\documentclass[onecolumn,11pt,draft=false]{IEEEtran}

\usepackage{enumitem}
\usepackage{balance}
\usepackage{setspace}
\usepackage{latexsym}

\usepackage{lipsum}
\usepackage{amsfonts}
\usepackage{epstopdf}
\usepackage{algorithmic}
\usepackage{algorithm}

\newtheorem{problem}{Problem}
\newtheorem{lemma}{Lemma}
\newtheorem{proposition}{Proposition}
\newtheorem{theorem}{Theorem}
\newtheorem{definition}{Definition}

\usepackage{amssymb,amsmath,bbm}
\usepackage{color}

\usepackage{cite}
\usepackage{amsmath,amssymb,amsfonts}
\usepackage{mathtools}
\usepackage{subfigure}
\usepackage{graphicx}
\usepackage{textcomp}
\usepackage{xcolor}
\def\BibTeX{{\rm B\kern-.05em{\sc i\kern-.025em b}\kern-.08em
    T\kern-.1667em\lower.7ex\hbox{E}\kern-.125emX}}

\newcommand{\realnonnegative}{{\mathbb{R}}_{\ge 0}}

\newcommand{\integerspositive}{\mathbb{Z}_{>0}}

\newcommand{\oprocendsymbol}{\hbox{$\bullet$}}
\newcommand{\oprocend}{\relax\ifmmode\else\unskip\hfill\fi\oprocendsymbol}

\DeclareMathOperator*{\argmin}{arg\,min}

\newcommand{\longthmtitle}[1]{\mbox{}\textup{\bf (#1):}}

\newcommand{\boldB}{\mathbf{B}}
\newcommand{\boldR}{\mathbf{R}}
\newcommand{\tI}{\tilde{I}}

\newcommand{\Nsum}{\sum_{i=1}^N}
\newcommand{\tsum}{\sum_{t=1}^n}

\newcommand{\tS}{\tilde{S}}
\newcommand{\tU}{\tilde{U}}

\newcommand{\bR}{\bar{R}}
\newcommand{\bB}{\bar{B}}
\newcommand{\bT}{\bar{T}}

\newcommand{\budget}{\mathcal{B}}
\newcommand{\budgetB}{\mathcal{B}_b}
\newcommand{\budgetR}{\mathcal{B}_r}
\newcommand{\F}{\mathcal{F}}
\newcommand{\N}{\mathcal{N}}
\newcommand{\G}{{\mathcal{G}}}
\newcommand{\E}{{\mathcal{E}}}
\newcommand{\X}{{\mathcal{X}}}
\newcommand{\Y}{{\mathcal{Y}}}

\newcommand{\bx}{\mathbf{x}}
\newcommand{\by}{\mathbf{y}}

\begin{document}

\title{\LARGE\bf Initialization and Curing Policies for Polya Contagion Networks{\LARGE $^*$}}

  \author{Greg Harrington,
         Fady Alajaji and
         Bahman Gharesifard}
  \renewcommand{\thefootnote}{\fnsymbol{footnote}}
  \footnotetext{$^*$Department of Mathematics and Statistics, Queen's University, Kingston, Ontario K7L 3N6, Canada (\texttt{12gth@queensu.ca}, \texttt{fa@queensu.ca}, \texttt{bahman.gharesifard@queensu.ca}). \\
  The work of the authors was funded by the National Science Foundation of Canada. The last author wishes to thank 
the Alexander von Humboldt Foundation for their generous support.}

\maketitle

\begin{abstract}
This paper investigates optimization policies for resource distribution in network epidemics using a model that derives from the classical Polya process. The basic mechanics of this model, called the Polya network contagion process, are based on a modified urn sampling scheme that accounts for both temporal and spatial contagion between neighbouring nodes in a network. We present various infection metrics and use them to develop two problems: one which takes place upon initialization and one which occurs continually as the Polya network process develops. We frame these problems as resource allocation problems with fixed budgets, and analyze a suite of potential policies. Due to the complexity of these problems, we introduce effective proxy measures for the average infection rate in each case. We also prove that the two-sided infection-curing game on the so-called expected network exposure admits a Nash equilibrium. In both the curing and initialization scenarios, we introduce heuristic policies that primarily function on the basis of limiting the number of targeted nodes within a particular network setup. Simulations are run for mid-to-large scale networks to compare performance of our heuristics to provably convergent gradient descent algorithms run on the simplified proxy measures.
\end{abstract}

\begin{IEEEkeywords}
Network epidemics; interacting Polya urn models; processes with reinforcement; temporal and spatial contagion; initialization and curing policies.
\end{IEEEkeywords}


\begin{spacing}{1.35}

\section{Introduction}

The study of epidemics on networks is an active research topic, e.g., see~\cite{DE-JK:10,PVM-JO-RK:09,CN-VMP-GJP:16} and the references therein. Real-life examples include the propagation of burst errors in a wireless communication channel~\cite{FA-TF:94}, of a biological disease through a population~\cite{LK-MA-KD-SK-AO:14}, of malware in computer or smartphone systems~\cite{MG-WG-DT:03}, and the dissemination of rumors~\cite{ER:03} and competing opinions~\cite{EA-LAA:05} in social networks.

Controlling the spread of contagion in networks has been thoroughly investigated for various systems, including the well-known Susceptible-Infected-Susceptible (SIS) and Susceptible-Infected-Recovered (SIR) epidemics models~\cite{PVM-JO-RK:09} and extensions~\cite{CN-VMP-GJP:17}, models based on real epidemic data~\cite{LW-JW:18}, as well as cascade models~\cite{DE-JK:10}.
Other problem formulations address link removal via immunization~\cite{NM-TK:04}, optimization of curing resources~\cite{CB-JC:10}, message-passing methods~\cite{MK-SS:11}, optimal control under performance and resource usage trade-offs~\cite{VP-MZ:14,KK-JK:14}, and competing or co-evolving contagions within networks~\cite{PEP-JL-CLB-AN-TB:17,WM-FB:17,QL-LZ:18,WW-QL:19,ARH-SS:19}.
We herein study initialization and curing strategies and related game-theoretic problems for the (infinite-memory) Polya network contagion model introduced in~\cite{MH-FA-BG:17, MH-FA-BG:18}.

The classical Polya process~\cite{GP-FE:23, GP-FE:28, GP:30} is a temporal contagion process that evolves through a method of sampling from a single urn containing a finite number of red and black balls, representing units of ``infection'' and ``healthiness,'' respectively. This model has been used in applications such as consensus dynamics~\cite{AF-AJ:11} and generalized for various other purposes~\cite{MC-MK:13,FC-SH-DJ:03}. The Polya network contagion process~\cite{MH-FA-BG:17, MH-FA-BG:18}, generalizes this sampling scheme to a general network, where each node is equipped with its own urn, by allowing for spatial interactions between neighbouring nodes. This is realized by introducing ``super urns," which are allocated for each network node and contain red and black balls from both a node's individual urn as well as those of neighbouring nodes~\cite{MH-FA-BG:17, MH-FA-BG:18}. In addition to temporal contagion, this modification generates spatial contagion not present in the classical model. 

The resulting Polya model for network contagion~\cite{MH-FA-BG:17, MH-FA-BG:18} is similar to the SIS model~\cite{DE-JK:10,YW-DC-CW-CF:03,HA-BH:13}. The Polya model can be operated in two modes: an infinite-memory mode where the reinforcing balls added at each time step remain permanently
in the underlying urn of each node, and a finite-memory mode where the added balls at a given time step remain in the system only for a fixed number of future time steps (this mode results in a finite-order Markov network contagion process).
While it was empirically observed in~\cite{MH-FA-BG:18} that the Polya network contagion model can mimic the behaviour of the SIS process, particularly in the finite-memory (Markovian) mode, there are some notable differences between the two models. 
For example, the Polya model benefits from exactly computable expressions for the joint and marginal probabilities of infection of its discrete-time contagion process, while the underlying (or exact) Markov process of the SIS model is typically represented and analyzed via a dynamical system using mean-field approximations; e.g., see \cite{PVM-JO-RK:09,HA-BH:13,NAR-HJA-BH:16,WM-SM-SZ-FB:17}. Furthermore, the Polya model captures the network contagion process at a microscopic level by associating an urn to each node, whose composition of red and black balls at each time step represents a granular profile for the degree of infection and healthiness of that node, which upon interaction with urns of nodes in proximity, results in infected/healthy states with each super urn draw. This allows us to minutely and analytically capture the stochastic evolution of the network-wide contagion process. Finally, the underlying discrete-time Markov chain of the SIS model has an absorbing (all healthy or disease-free) state~\cite{NAR-HJA-BH:16}, while the Polya model does not in general\footnote{The Polya model is more general than the exact SIS Markov process in terms of allowing different asymptotic behaviours for the contagion process. For example, in the finite-memory $M$ mode with homogeneous parameters, the Polya network contagion process is an $M$'th order irreducible and aperiodic Markov chain with a unique stationary distribution that is not necessarily the all-healthy state with probability one.} as its contagion process is symbiotically generated from its initial state (of red/black ball mixtures in the individual urns of the network nodes) via the Polya ball sampling mechanism.


In this paper, we examine resource allocation problems for the (infinite-memory) Polya network contagion process.
The first problem considers a one-time application of a control policy upon initialization of the network, while the second concerns the continual distribution of resources as the network contagion process develops. 
More specifically, the initialization problem concerns one party judiciously controlling the distribution of resources in an effort to minimize the average network-wide infection rate. This is a novel {\em preemptive} infection mitigation problem, not typically considered in SIS models, due to the Polya model's inherent characteristics where contagion organically develops and propagates from initial network conditions.
Resources in this case concern the distribution of red and black balls within a network, prior to any draws taking place. We focus on the one-sided finite horizon case, wherein a player controls the allocation of black balls according to a fixed budget, with the goal of minimizing the average infection rate at some fixed point in the future. We provide a series of results showing that optimal policies for this problem satisfy two conditions: nested nodes will receive no resources, and symmetric red ball initializations will yield symmetric black ball initializations (the specifics of what is meant by ``symmetry" in this case is further detailed using fundamental ideas from graph theory). Since determining the average infection rate at time $n\ge 1$ becomes increasingly complex as $n$ increases, we use the average infection rate at time $n=1$ to simplify the problem, and we prove that under certain conditions this proxy measure admits an optimal policy. An algorithm to obtain such an optimal policy using gradient descent is given and used as a base with which we compare our own heuristic policies. We detail three varieties of heuristic policies, each of which works on the basis of limiting our set of viable target nodes. The first set of policies target the set of ``inner nodes", the second works through a layered application of this inner node targeting technique, and the last works by iteratively targeting the most central nodes until full network coverage (either direct or indirect) is obtained. We provide two different algorithms for determining these target sets and variations based on previously tested heuristics (see~\cite{MH-FA-BG:19}).

The second control problem, referred to as the ``infection-curing problem" or the ``Delta-curing problem," concerns {\em reactive} intervention policies deployed during the progress of contagion in the network. In this setup, which is an extension of the contagion curing problem studied in~\cite{MH-FA-BG:19},  we define a two-player game on the average infection rate for a given network. Prior work on competitive dynamics over SIS and social networks include~\cite{PEP-JL-CLB-AN-TB:17} and~\cite{WM-FB:17}, respectively.  
In order to simplify the problem of finding optimal control policies, we consider a game on the proxy measure of the expected network exposure. We show that the expected network exposure is convex as a function of the curing parameters and concave as a function of the infection parameters. We prove that under budget constraints, there exists a deterministic Nash equilibrium~\cite{JN:50} for the game on the expected network exposure, which can be determined using gradient descent algorithms. We also develop a set of heuristic policies for the Delta-curing problem using the same methods as in the initialization case.
We next provide simulations to demonstrate the performance of the proposed heuristics for the initialization and curing problems on both sparse and dense large-scale networks. These policies are compared to optimal policies for the one-step proxy measures. 

The rest of the paper is organized as follows. In Section~\ref{ch:network_process}, we describe in sufficient detail
the (infinite-memory) Polya network contagion process~\cite{MH-FA-BG:18}.
The initialization and infection-curing problems are theoretically investigated in Sections~\ref{ch:initialization} and~\ref{ch:curing-infection}, respectively. Heuristic optimization strategies for the above problems are developed in Section~\ref{sec:heuristics} and simulation results are presented in
Section~\ref{ch:simulations}. Finally, conclusions are stated in ~Section~\ref{ch:conclusion}.



\section{Polya Network Contagion}\label{ch:network_process}
We devote this section to some background on the classical Polya urn process and the (infinite-memory) Polya network contagion process under study. Throughout, a number of results are derived using standard probability concepts, which can be found in various texts including~\cite{RA-CD:00,GG-DS:01}. Herein, let $(\Omega,\F,P)$ be a probability space. Consider the stochastic process $\{Z_n\}_{n=1}^{\infty}$, where $Z_n=(Z_{1,n},\ldots,Z_{N,n})$ is a random vector on $\Omega$. Here we think of the process indices as time indices, while the vector indices delineate spatial indices (to distinguish nodes within a network). For $i=1,\ldots,N$, we write the $n$-tuple $(Z_{i,1},\ldots,Z_{i,n})$ as $Z_i^n$, and more generally, we write $Z_{i,s}^t=(Z_{i,s}, Z_{i,s+1},\ldots,Z_{i,t})$ for $s<t$.

\subsection{Classical Polya Process~\cite{GP-FE:23, GP-FE:28, GP:30}}

An urn initially contains $T=R+B$ balls of which $R \in\integerspositive$ are red and $B \in \integerspositive$ black. At each time, $n$, a ball is drawn from the urn and returned with $\Delta>0$ balls of the same colour. The random variable $Z_{n}$ is used to indicate the colour of the ball on the $n$th draw, where
\begin{align*}
	Z_n = \begin{cases}
		1 &\text{if the $n$th draw is red}\\
		0 &\text{if the $n$th draw is black.}
	\end{cases}
\end{align*}
We can then use $U_{n}$ to denote the proportion of red balls in the urn after the $n$th draw, where
\begin{align*}
	U_n &:= \frac{R + \Delta\tsum Z_{t}}{T + n\Delta} = \frac{\rho_c + \delta_c\tsum Z_{t}}{1+n\delta_c}.
\end{align*}
Here $\rho_c = \frac{R}{T}$ is the initial proportion of red balls in the urn and $\delta_c = \frac{\Delta}{T}$ is a correlation parameter. The conditional probability of drawing a red ball at time $n$ , given the past history of draws $Z^{n-1} := (Z_{1}, \ldots, Z_{n-1})$, is given by
\begin{align}\label{eq:U_n}
	P(Z_n = 1 \ | \ Z^{n-1} ) &= \frac{R + \Delta\sum_{t=1}^{n-1} Z_{t}}{T + (n-1)\Delta} = U_{n-1}.
\end{align}
Additionally, we note that the  process $\{Z_n\}_{n=1}^{\infty}$ is exchangeable (thus stationary) and non-ergodic and that
both $U_n$ and the process sample average $\frac{1}{n}\sum_{t=1}^nZ_t$ converge almost surely (a.s.) to a random variable governed by a Beta distribution with parameters $\frac{\rho_c}{\delta_c}$ and $\frac{1-\rho_c}{\delta_c}$~\cite{GP:30,WF:71}.
Variations on this model have been explored in a number of other works. We focus exclusively on an adaptation of this model first introduced in~\cite{MH-FA-BG:17} and further explored in~\cite{MH-FA-BG:18},\cite{MH-FA-BG:19}, designed as a network process that introduces spatial contagion to the classical temporal contagion process presented above.

\subsection{Polya Network Contagion Process}
We now review the Polya network contagion process~\cite{MH-FA-BG:17, MH-FA-BG:18}. 
Let $\G = (V,\E)$ be a graph, where $V = \{1, \ldots, N\}$ is the set of $N \in \integerspositive$ nodes and $\E \subset V\times V$ is the set of edges. Throughout, $\G$ is taken to be undirected, i.e., $(i,j)\in\E$ if and only if $(j,i)\in\E$ for all $i,j\in V$, and is assumed to be connected, i.e., there exists a series of edges connecting any two nodes in $\G$. The set of neighbours to node $i$ is given by $\N_i=\{v\in V:(i,v)\in\E\}$. We then define $\N_i'=\{i\}\cup\N_i$. Each node $i\in V$ is assigned an urn with $R_i \in  \mathbb{Z}_{\ge 0}$ red balls and $B_i \in \mathbb{Z}_{\ge 0}$ black balls, with the total number of balls at node $i$ denoted by $T_i = R_i+B_i$. We denote the black and red ball initializations by the vectors $\boldB:=(B_1,\dots, B_N)$ and $\boldR:=(R_1,\dots,R_N)$, respectively. We will refer to $\boldB$ as the curing initialization and $\boldR$ as the infection initialization.

\begin{figure}[!ht]
\centering 
	\includegraphics[width=0.45\linewidth]{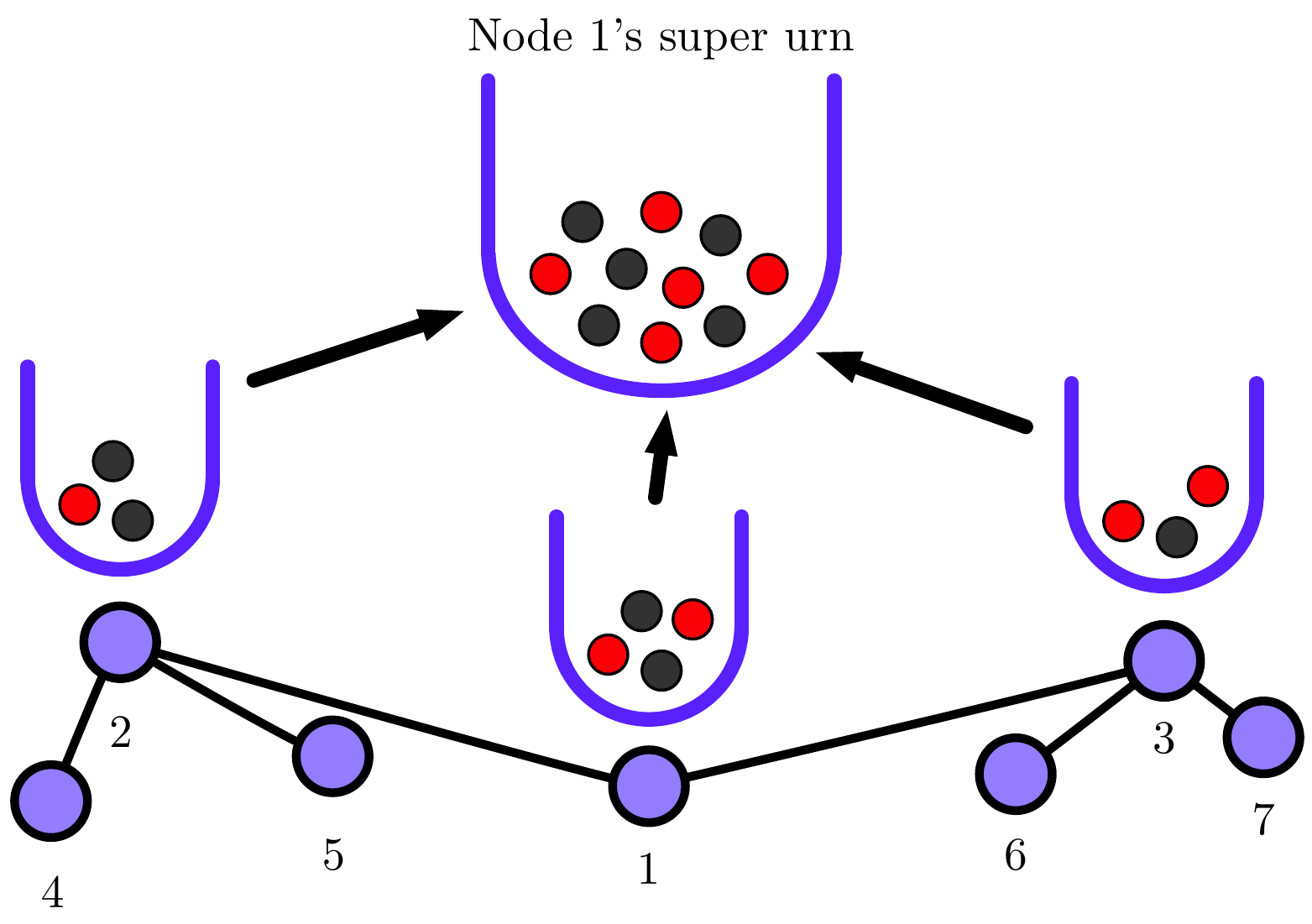}
	\caption[Illustration of a super urn in a network]{Illustration of a super urn in a network~\cite{MH:17}.}
	\label{fig:super_urn}
\end{figure}

We next introduce the concept of a super urn for a given node $i$, which contains the red and black balls at node $i$ as well as the red and black balls of each of its neighbouring nodes. This concept is illustrated in Figure~\ref{fig:super_urn}. For each node $i\in V$, we use $\bR_i=\sum_{j \in \N_i'}R_j$ and $\bB_i=\sum_{j \in \N_i'}B_j$ to denote the number of red and black balls, respectively, in node $i$'s super urn, and we assume throughout that $\bR_i + \bB_i > 0$ for all $i \in V$ (to ensure that the super urn associated to each node is non-empty). The total number of balls in the $i$th super urn is given by $\bT_i = \bR_i+\bB_i$. Similar to the classical Polya process, a draw is conducted for each node at each time, and then a number of balls of the same colour is added to that node's individual urn. In this case, however, the draw is conducted on the super urn of each node. As well, we may allow for the number of added balls to vary based on which node the draw was for, the colour of the draw, and the time at which the draw occurred. Thus, for node $i$ at time $n$, if a red ball is drawn, we add $\Delta_{r,i}(n)$ red balls to node $i$'s individual urn, and if a black ball is drawn, we add $\Delta_{b,i}(n)$ black balls to node $i$'s individual urn. We also assume that there exists $i \in V$ and $n$ such that $\Delta_{b,i}(n)+\Delta_{r,i}(n) > 0$. 
We let $Z_{i,n}$ indicate the colour of the $n$th draw for node $i$:
\begin{align*}
	Z_{i,n} = \begin{cases}
		1 &\text{if the $n$th draw for node $i$ is red}\\
		0 &\text{if the $n$th draw for node $i$ is black.}
	\end{cases}
\end{align*}
Translating the above parameters to the framework of an epidemic in a population, 
black and red balls represent units of ``healthiness''
(such as white blood cells) and ``infection'' (such as bacteria or viruses), respectively~\cite{MH-FA-BG:18}. In the super urn of a given node (such as a ``person'' or an ``organism''), white blood cells (resp., bacteria) combine 
to improve (resp., impair) the overall health in the node's neighborhood. Drawing red at time $n$ from the super urn means that the bacteria in the neighborhood were successful in replicating, making the individual more infected or less immune; 
otherwise, they were healthier since the white blood cells reproduced (as balls of the the same color as the drawn balls are added). Thus, when $Z_{i,n} = 1$, we declare that node $i$ is infected at time $n$, and if $Z_{i,n} = 0$, then it is healthy.
We refer to  $\{\Delta_{b,i}(n)\}_{n=1}^{\infty}$ as the curing parameters and  $\{\Delta_{r,i}(n)\}_{n=1}^{\infty}$ as the infection parameters for node $i$.
From an urn-sampling point of view, these $\Delta$'s are nonnegative integers; however, we allow them to be nonnegative real-valued numbers for mathematical convenience (e.g., when determining optimal control policies). We will likewise do the same for the initialization parameters $\boldB$ and $\boldR$.
Ultimately, the collection of draw values at each time step form our Polya network contagion process, represented as $\{Z_n\}_{n=1}^{\infty}$, where $Z_n$ represents the vector $(Z_{1,n},\dots,Z_{N,n})$. If desired, we can also separate out the individual draw process $\{Z_{i,n}\}_{n=1}^{\infty}$ for a given node $i$. 
We introduce some metrics to measure the infection spread in the network.

Similar to the classical Polya process, we use $U_{i,n}$ to denote the proportion of red balls in node $i$'s urn after the $n$th draw. Letting
\begin{align}\label{eq:X_i,n}
X_{i,n}=T_i + \sum_{t=1}^{n} \Delta_{r,i}(t)Z_{i,t} + \Delta_{b,i}(t)(1-Z_{i,t})
\end{align}
with $ T_i>0 $, represent the total number of balls in node $i$'s urn after the $n$th draw, we can write
\begin{align}\label{eq:U_i,n}
	U_{i,n} &= \frac{R_i + \tsum \Delta_{r,i}(t)Z_{i,t}}{X_{i,n}}.
\end{align}
We note that $U_{i,0} = \frac{R_i}{T_i}$ represents the initial (at time $n=0$) proportion of red balls in node $i$'s individual urn. Though for the classical Polya process the individual urn proportion is commonly used in the analysis, in this case, the information we gain through this value is limited by the fact that our draws are performed on the super urns for each node, rather than the urns for the individual nodes themselves. We therefore also need a metric to measure the proportion of red balls in the super urn of node $i$ at time $n$, which we denote by $S_{i,n}$. This proportion is given by the equation
\begin{align}\label{eq:S_i,n}
	S_{i,n} 
	&= \frac{\bR_i + \sum_{j \in \N_i'}\sum_{t=1}^{n} \Delta_{r,j}(t)Z_{j,t}}{\sum_{j \in \N_i'}X_{j,n}}=\frac{\sum_{j \in \N_i'}U_{j,n}X_{j,n}}{\sum_{j \in \N_i'}X_{j,n}}.
\end{align}
Like for the individual node proportions at time $n=0$, we use the convention $S_{i,0}=\frac{\bR_i}{\bT_i}$. This is simply representative of the proportion of red balls in node $i$'s super urn when the process is initialized. We next define
\begin{align*}
	\tU_n &= \frac{1}{N}\sum_{i=1}^N U_{i,n},
\end{align*}
which we refer to as the \textit{network susceptibility}, and
\begin{align*}
	\tS_n &= \frac{1}{N}\sum_{i=1}^N S_{i,n},
\end{align*}
which we refer to as the \textit{network exposure}. Both the network susceptibility and network exposure are functions of the underlying network contagion process $\{Z_n\}_{n=1}^{\infty}$, as are the respective measures for each node. In general, for notational ease, we will omit these arguments, unless otherwise noted. For example, for a given node $i$, we will write $S_{i,n}$ rather than $S_{i,n}(\{Z_j^n\}_{j=1}^N)$.

We next calculate infection probabilities within the network, i.e., the likelihood of drawing red balls. To start, the conditional probability of drawing a red ball for node $i$ at time $n$ , given the past history of draws for all the nodes 
\begin{equation*}
	\{Z_j^{n-1}\}_{j=1}^N := \{(Z_{1,1}, \ldots, Z_{1,n-1}),\ldots, (Z_{N,1}, \ldots, Z_{N,n-1})\},
\end{equation*}
is given by
\begin{align}\label{eq:conditional}
	P\left(Z_{i,n} = 1 | \{Z_j^{n-1}\}_{j=1}^N\right) &= \frac{\bR_i + \sum_{j \in \N_i'}\sum_{t=1}^{n-1} \Delta_{r,j}(t)Z_{j,t}}{\sum_{j \in \N_i'}X_{j,n-1}} = S_{i,n-1}.
\end{align}
This is analogous to~\eqref{eq:U_n} of the classical Polya process, except that for the network contagion process this conditional probability is represented by the super urn proportion rather than the individual urn proportion. We also need an unconditioned version of this equation. The $n$-fold joint probability of the network $\G$, which is a function defined for each specific draw history $a_i^n\in\{0,1\}^n$, $i\in V$, can be written as
\begin{align}\label{eq:joint_prob}
	P_{\G}^{(n)}(a_1^n,\dots,a_N^n)&:=P(\{Z_i^n=a_i^n\}_{i=1}^N)\nonumber\\
									&=\prod_{t=1}^nP(\{Z_{i,t}=a_{i,t}\}_{i=1}^N|\{Z_i^{t-1}=a_i^{t-1}\}_{i=1}^N)\nonumber\\
									&=\prod_{t=1}^n\prod_{i=1}^N(S_{i,t-1})^{a_{i,t}}(1-S_{i,t-1})^{1-a_{i,t}}.
\end{align}
It is important to note here that the draw value of different nodes at a given time $t$ are independent when conditioned on the history of the process up to that point, thus allowing us to write $P_{\G}^{(n)}$ in the form above. Finally, we introduce the metric
\begin{align}\label{eq:infection_rate}
	\tI_n:=\frac{1}{N}\Nsum P(Z_{i,n}=1),
\end{align}
which we call the \textit{average infection rate at time n}. Since $P(Z_{i,n}=1)$ represents the marginal probability that a node $i$ is infected at time $n$, the average infection rate gives us the network-wide average of seeing infection within a node. Thus, this value represents the marginal probability that a randomly selected node within the network is infected at time $n$. In most cases where we look to study the asymptotic behaviour of the network contagion process, we do so in regard to this metric.

We can further break down~\eqref{eq:infection_rate} by noting that
\begin{align}\label{eq:marginal_node}
	P(Z_{i,n}=1) &= \sum_{\{a_j^{n-1}\}_{j=1}^N}P(Z_{i,n}=1|\{Z_j^{n-1}=a_j^{n-1}\}_{j=1}^N)P(\{Z_j^{n-1}=a_j^{n-1}\}_{j=1}^N) \nonumber\\
				&= \sum_{\{a_j^{n-1}\}_{j=1}^N}S_{i,n-1}(\{a_j^{n-1}\}_{j=1}^N)P_{\G}^{(n-1)}(\{a_j^{n-1}\}_{j=1}^N),
\end{align}
where the summation is over all possible draw histories $a_j^{n-1}\in\{0,1\}^{n-1}$, $j\in V$.

Each of these measures, $\tU_n$, $\tS_n$, and $\tI_n$,are closely related to one another and can be useful in various situations. However, all of these values are network-specific as in general all three depend in some way on the underlying network topology and the initial ball distributions. The quantity $\tI_n$ is particularly difficult to analyze, and hence we often resort to using alternative measures that allow us to simplify our analysis or gain insight into how $\tI_n$ behaves.

By inspecting~\eqref{eq:U_i,n} and~\eqref{eq:S_i,n}, it is clear that an increase in the node-specific proportion of red balls increases the super urn proportion of red balls in all neighbouring nodes. Likewise, if we consider how~\eqref{eq:conditional} relates to~\eqref{eq:marginal_node}, it would make sense that an increase in the super urn proportion would also lead to an increase in the marginal network-wide infection rate; we further discuss these observations in later sections. For now, we have that
\begin{align*}
	\uparrow U_{i,n}\implies \uparrow S_{j,n}\text{ for all }j\in \N_i',
\end{align*}
and since $\tI_n$ is essentially the marginal version of the conditional $\tS_n$, we often start our analysis by using $\tS_n$ in place of $\tI_n$, without losing significant information.

\section{Initialization Problems}\label{ch:initialization}

When we consider the problem of limiting the spread of infection for the Polya network contagion model, there are different ways to implement a curing policy. One method is to control the allocation of our curing parameters $\{\Delta_{b,i}(n)\}_{i=1}^N$ at each time $n$, which we explore in Section~\ref{ch:curing-infection}. We can, alternatively, control how the process is initialized. This is where we currently draw our focus for discussion.

The initialization parameters $\boldB\in\realnonnegative^N$ and $\boldR\in\realnonnegative^N$ can be tailored to alter the evolution of the Polya network contagion process for a given network. The goal is to determine how to allocate such resources so as to achieve a more desirable result, and determine whether or not a theoretically optimal initial allocation policy exists. Such optimizations may be performed for either or both sets of initialization parameters subject to an initialization budget $\budget\in\realnonnegative$. We start with a formal setup for the one-sided initialization problem, wherein we aim to minimize the average infection rate for a given time $n$ over the black ball initialization parameters $\boldB$ subject to a budget $\budgetB$. Throughout this section we assume the initial distribution of red balls $\boldR$ is fixed
such that $\bR_i > 0$ for all $i \in V$.

\begin{problem}\longthmtitle{One-Sided Finite Horizon Constrained Budget Initialization}\label{pb:init1}
For a fixed time $n$, minimize the average infection rate $\tI_n$ subject to a budget $\budgetB$ on the curing initialization $\boldB=(B_1,\dots,B_N)$, i.e., find:
\begin{align*}
	\min_{\{B_i\}_{i=1}^N\in\realnonnegative^N:\sum_{i=1}^NB_i\leq\budgetB}\tI_n.
\end{align*}
\end{problem}

We remark that the problem above can be extended to the infinite horizon case, but for the purposes of our work, we limit ourselves to a finite horizon. Moreover, one can establish a \emph{two-sided} initialization problem, wherein the objective is to minimize the average infection rate for a given $n$, but we do so over both sets of initialization parameters $\boldB$ and $\boldR$. We, however, only focus on Problem~\ref{pb:init1} in this work.
We first consider Problem~\ref{pb:init1} for $n=1$. We have
\begin{align}\label{eq:infect_1}
	\tI_1	&=\frac{1}{N}\sum_{i=1}^NP(Z_{i,1}=1) = \frac{1}{N}\sum_{i=1}^N\frac{\bR_i}{\bR_i+\bB_i}.
\end{align}
As a function of the initialization parameters, $\tI_1$ is well-defined over $\{(\boldB,\boldR)\in\realnonnegative^{2N}|\bB_i+\bR_i\neq0, \forall i\in V\}$. Furthermore, $\tI_1$ is convex in $\boldB$ and concave in $\boldR$.
We use the form of~\eqref{eq:infect_1} to provide some simplifications regarding solutions to Problem~\ref{pb:init1} for $n=1$. We will then extend these results to a general time $n$.

\subsection{Outer Node Allocation} 
We first show that for the one-sided initialization problem (Problem~\ref{pb:init1}) with $n=1$, nodes with a nested neighbourhood (i.e., nodes with neighborhood $\N'$ contained in the neighborhood of another node) can be ignored. We refer to such nodes as ``outer nodes" due to their topological location.

\begin{lemma} \label{lma:313}
For a general network $\G=(V,\E)$ equipped with the Polya network contagion model, if for any nodes $i,j\in V$ we have $\N_i' \subset \N_j'$, then it can be assumed that $B_{i}=0$ when minimizing $\tI_1$ over the initial distribution of black balls. Likewise, we can assume $R_i=0$ when maximizing $\tI_1$ over the initial distribution of red balls.
\end{lemma}

\begin{IEEEproof}Consider a given initial distribution of black balls $\mathbf{B}=(B_1,\ldots,B_N)$, and red balls $\mathbf{R}=(R_1,\ldots,R_N)$. As written in~\eqref{eq:infect_1}, we have
\begin{align*}
	\tI_1 = \frac{1}{N}\sum_{k=1}^N\frac{\bR_k}{\bR_k+\bB_k}.
\end{align*}
Without loss of generality, assume that $\N_1'\subset\N_2'$, and let the distribution  $\mathbf{B}^*=(0,B_1+B_2,B_3,\ldots,B_N)$ be an alternative curing initialization. For a given node $k\in V$, we have 3 cases to consider:
\begin{enumerate}
	\item $k\in\N_1'\implies k\in\N_2'$;
	\item $k\notin\N_1'$ and $k\in\N_2'$;
	\item $k\notin\N_2'\implies k\notin\N_1'$.
\end{enumerate}
For cases 1 and 3, note that $\sum_{l\in\N_k'}B_l=\sum_{l\in\N_k'}B_l^*$. For case 2, we have $\sum_{l\in\N_k'}B_l\leq\sum_{l\in\N_k'}B_l^*$ with equality if and only if $B_1=0$. Denoting $\bB_k^*=\sum_{l\in\N_k'}B_l^*$, we have
\begin{align*}
	\frac{\bR_k}{\bR_k+\bB_k}\geq\frac{\bR_k}{\bR_k+\bB_k^*}\text{, $\forall k\in V$}.
\end{align*}
Therefore, for $\tI_1$ as a function of the initial distribution of black balls we have that $\tI_1(\mathbf{B}^*)\leq\tI_1(\mathbf{B})$. The proof for the maximization case follows similarly.
\end{IEEEproof}

To extend Lemma~\ref{lma:313} for general $n$, we derive the following intermediate results.

\begin{lemma}\label{lma:314}
Consider two sequences of draw values $\{a_{j}^{n}\}_{j=1}^{N}$ and $\{b_{j}^{n}\}_{j=1}^{N}$, where both sequences are equal, except for some $(k,s)\in V\times\{1,\ldots,n\}$, where $a_{k,s}=1$ and $b_{k,s}=0$. Then $P(Z_{i,n+1}=1|\{Z_{j}^{n}=a_{j}^{n}\}_{j=1}^{N})\geq P(Z_{i,n+1}=1|\{Z_{j}^{n}=b_{j}^{n}\}_{j=1}^{N})$, for all $i\in V$.
\end{lemma}

\begin{IEEEproof}
Assume that $k\in\mathcal{N}_{i}'$ (note that if this is not the case, the two values we wish to compare are equal). Let
\begin{align*}
	y=\sum_{j\in\mathcal{N}_{i}'}\sum_{t=1}^{n}a_{j,t}\Delta_{r,j}(t)\ \quad
	x=\sum_{j\in\mathcal{N}_{i}'}\sum_{t=1}^{n}(1-a_{j,t})\Delta_{b,j}(t).
\end{align*}
Then letting $y^*=y-\Delta_{r,k}(s)$, and $x^*=x+\Delta_{b,k}(s)$, we obtain
\begin{align*}
	P(Z_{i,n+1}=1|\{Z_{j}^{n}=a_{j}^{n}\}_{j=1}^{N})
		&\geq \frac{\bR_i+y}{\bR_i+\bar{B}_{i}+y+x^*}\\
		&\geq \frac{\bar{R}_{i}+y^*}{\bar{R}_{i}+\bar{B}_{i}+y^*+x^*}\\
		&=P(Z_{i,n+1}=1|\{Z_{j}^{n}=b_{j}^{n}\}_{j=1}^{N}).
\end{align*}
\end{IEEEproof}

Intuitively, this result states that the process is self-reinforcing; i.e., higher rates of infection increase the likelihood of the infection recurring. We will use this lemma in the proof of Lemma~\ref{lma:315}.

Given the general form of $\tI_n$ given in~\eqref{eq:infection_rate}, we next compare the performance of two curing initializations, $\boldB$ and $\boldB^*$, in limiting the spread of infection. In order to differentiate between the resultant probabilities, we will use $^*$ to denote any probability obtained using the curing initialization $\boldB^*$, e.g., $\tI_n^*=\tI_n(\boldB^*)$.


\begin{lemma}\label{lma:315}
For a general network $\G=(V,\E)$ equipped with the Polya network contagion model with two different curing initializations $\mathbf{B}=(B_1,\dots,B_N)$ and $\mathbf{B}^*=(B_1^*,\dots,B_N^*)$. If $S^*_{i,t}(a_1^t,\ldots,a_N^t)\leq S_{i,t}(a_1^t,\ldots,a_N^t)$ for all $(i,t)\in V\times\{1,\ldots,n\}$ and any sequence of draws $\textbf{a}^n=(a_{1,s},\dots,a_{N,s})_{s=1}^n\in\{0,1\}^{N\times n}$, then $\tI_{n+1}^*\leq\tI_{n+1}$.
\end{lemma}

\begin{IEEEproof}
Let $\{Y_{i,t}: t=1,\dots,n,i=1\dots,N\}$ be independent random variables each uniformly distributed on the unit interval [0,1]. For all $i=1,\dots,N$, let $Z_{i,t}$ and $Z_{i,t}^*$ be inductively defined for $t=1,\dots,n$ by
\begin{equation}\label{eq:Z_alternate}
	Z_{i,t} = \begin{cases}
		1 &\text{if $Y_{i,t}\leq S_{i,t-1}(\{Z_i^{t-1}\}_{i=1}^N)$}\\
		0 &\text{otherwise,}
	\end{cases}
\end{equation}
and
\begin{equation}\label{eq:Z*_alternate}
	Z_{i,t} = \begin{cases}
		1 &\text{if $Y_{i,t}\leq S_{i,t-1}^*(\{Z_i^{*,t-1}\}_{i=1}^N)$}\\
		0 &\text{otherwise.}
	\end{cases}
\end{equation}
Here we assume that $S_{i,0}\geq S_{i,0}^*$ for all $i=1\dots,N$.
Note that since 
\begin{align*}
&P(Z_{i,t}=1|\{Z_i^{t-1}\}_{i=1}^N)=S_{i,t-1}(\{Z_i^{t-1}\}_{i=1}^N), \\ 
\mathrm{and} \quad &P(Z_{i,t}^*=1|\{Z_i^{*,t-1}\}_{i=1}^N)=S_{i,t-1}^*(\{Z_i^{*,t-1}\}_{i=1}^N),
\end{align*}
for all $i=1,\dots,N$ and $t=1,\dots,n$, the processes have the required Polya network process distribution with initializations $\boldB$ and $\boldB^*$, respectively.

We claim that if for some $t\in\{1,\dots,n-1\}$ the sequence $\{Z_i^{t-1}\}_{i=1}^N$ dominates $\{Z_i^{*,t-1}\}_{i=1}^N$, in the sense that $Z_{i,s}\geq Z_{i,s}^*$ for all $i=1,\dots,N$ and $s=1,\dots,t-1$, then $\{Z_i^{t}\}_{i=1}^N$ dominates $\{Z_i^{*,t}\}_{i=1}^N$. To prove this claim, note that by repeated application of Lemma~\ref{lma:314}, we have that
 \[
 S_{i,t-1}(\{Z_i^{t-1}\}_{i=1}^N)\geq S_{i,t-1}^*(\{Z_i^{*,t-1}\}_{i=1}^N),
 \] 
 for all $i=1,\dots,N$; specifically, since $Z_{i,s}\geq Z_{i,s}^*$ means that either $Z_{i,s}=Z_{i,s}^*$ or $Z_{i,s}=1$ and $Z_{i,s}^*=0$, Lemma~\ref{lma:314} is applied $K=|\{(i,s):Z_{i,s}\neq Z_{i,s}^*\}|$ times. This,~\eqref{eq:Z_alternate}, and~\eqref{eq:Z*_alternate} then give $Z_{i,t}\geq Z_{i,t}^*$, which in turn implies that $\{Z_i^{t}\}_{i=1}^N$ dominates $\{Z_i^{*,t}\}_{i=1}^N$, finishing the proof of the claim. 
 
 To proceed with the proof of the result, by using the conclusion of the claim for $t=0$, which vacuously holds by the assumption that $S_{i,0}\geq S_{i,0}^*$ for all $i=1,\dots,N$, we obtain $Z_{i,1}\geq Z_{i,1}^*$, $i=1,\dots,N$. Applying now the result of the claim inductively, for $t=2,\dots,n$, we obtain that
\begin{align*}
	Z_{i,t}\geq Z_{i,t}^*,\quad i=1,\dots,N, \quad t=1,\dots,n
\end{align*}
with probability one. Thus the event $\{Z_{i,t}^*=1\}$ implies $\{Z_{i,t}=1\}$ and so $P(Z_{i,t}^*=1)\leq P(Z_{i,t}=1)$, for all $i=1,\dots,N$ and $t=1,\dots,n$. Therefore $
	\tI_t^*\leq\tI_t^* $,
for all $t=1,\dots,n$.
\end{IEEEproof}

We can now use this result to extend Lemma~\ref{lma:313} to the case of a general time $n$.

\begin{theorem}\label{lma:316}
For a general network $\G=(V,\E)$ equipped with the Polya network contagion model, if for any nodes $i,j\in V$ we have $\N_i' \subset \N_j'$, then it can be assumed that $B_{i}=0$ when minimizing $\tI_n$ over the curing initialization. Likewise, it can be assumed that $R_i=0$ when maximizing $\tI_n$ over the infection initialization.
\end{theorem}

\begin{IEEEproof}
Without loss of generality, assume $\N_1' \subset \N_2'$. Consider the general form of $S_{i,n}(a_1^n,\ldots,a_N^n)$:
\begin{align}\label{eq:Sin}
	\frac{\bR_i+\sum_{j\in\N_i'}\sum_{t=1}^n a_{j,t}\Delta_{r,j}(t)}{\bR_i+\bB_i+\sum_{j\in\N_i'}\sum_{t=1}^n (a_{j,t}\Delta_{r,j}(t)+(1-a_{j,t})\Delta_{b,j}(t))}.
\end{align}
Letting $B_1^*=0$, and $B_2^*=B_1+B_2$, as in the proof of Lemma~\ref{lma:313}, we have for any $i\in V$ that $\bB_i^*\geq\bB_i$. Plugging into~\eqref{eq:Sin} yields that $S_{i,n}^*\leq S_{i,n}$ for any realization of the draw process. This will hold for any $i\in V$ and for any time $t$. Applying Lemma~\ref{lma:315} yields the desired result.
\end{IEEEproof}

%

\subsection{Symmetry}
The next simplification we make regarding solutions to Problem \ref{pb:init1} has to do with networks that possess ``symmetry'' properties, in a sense that we make precise shortly. We start with some preliminaries, referring the reader to~\cite{JAB-USRM:76, NB:94, CDG-GFR:01} for more details. 
\begin{definition}
An automorphism of a graph $\G=(V,\E)$ is a permutation $\sigma$ of $V$, such that $(u,v)\in\E$ if and only if $(\sigma(u),\sigma(v))\in\E$. The set of automorphisms of a given graph forms a group under the composition operation, called the automorphism group of $\G$ and denoted $Aut(\G)$.
\end{definition}
\begin{definition}
For a graph $\G=(V,\E)$, if $H\leq Aut(\G)$ is a subgroup of automorphisms of $\G$, then $u,v\in V$ are similar under $H$ if there exists an automorphism in $H$ which maps $u$ to $v$. Equivalence classes defined by similarity under $H$ are called orbits of the graph $\G$ by $H$. The sets of orbits by $H$ form a partition of the vertices of $\G$, called an orbit partition.
\end{definition}
If there exists an automorphism which maps between two different nodes, it means that those nodes are distinguishable only by their labels. 
For our purposes, ``symmetry" between nodes is characterized by graph automorphisms, and an orbit partition allows us to divide a network based on these symmetries. 
Our goal is to prove that under certain conditions, sets of nodes that are similar under some $H\leq Aut(\G)$ will be treated equivalently; i.e., they will all be given the same number of resources by an optimal policy. We start with a preliminary result.

\begin{lemma}\label{lma:lma2.5}
Consider a general network $\G=(V,\E)$ with a given initialization $(\boldB,\boldR)=(B_1,\dots, B_N,R_1,\dots,R_N)$. Suppose there exists a subset of nodes $V'\subset V$ such that $\bR_i=\bR_j$ for all $i,j\in V'$. Now suppose we can find some alternate curing initialization $\mathbf{B}^*=(B^*_1,\dots,B^*_N)$ such that $\bB_i^*=\frac{1}{|V'|}\sum_{j\in V'}\bB_j$ for all $i\in V'$. Then $\sum_{i\in V'}P^*(Z_{i,1}=1)\leq\sum_{i\in V'}P(Z_{i,1}=1)$.
\end{lemma}

\begin{IEEEproof}
We prove this result by induction on $|V'|$. Start with the base case, where $|V'|=2$. Consider two nodes $i,j\in V$ such that $\bR_i=\bR_j$. Let $\bB^*_i=\bB^*_j=\frac{1}{2}(\bB_i+\bB_j)$. We have
\begin{align}
	P(Z_{i,1}=1)+P(Z_{j,1}=1)	&=\frac{\bR_i}{\bR_i+\bB_i}+\frac{\bR_i}{\bR_i+\bB_j},\label{eq:eq5}\\
	P^*(Z_{i,1}=1)+P^*(Z_{j,1}=1)	&=\frac{2\bR_i}{\bR_i+\frac{1}{2}(\bB_i+\bB_j)}.\label{eq:eq6}
\end{align}
The difference between~\eqref{eq:eq5} and~\eqref{eq:eq6} yields
\begin{align*}
	\frac{\bR_i}{\bR_i+\bB_i}&+\frac{\bR_i}{\bR_i+\bB_j}-\frac{2\bR_i}{\bR_i+\frac{1}{2}(\bB_i+\bB_j)}\\
	=&\bR_i\left[\frac{2\bR_i+\bB_i+\bB_j}{(\bR_i+\bB_i)(\bR_i+\bB_j)}-\frac{2}{\bR_i+\frac{1}{2}(\bB_i+\bB_j)}\right]\\
	=&\bR_i\left[\frac{(2\bR_i+\bB_i+\bB_j)(\bR_i+\frac{1}{2}(\bB_i+\bB_j))-2(\bR_i+\bB_i)(\bR_i+\bB_j)}{(\bR_i+\bB_i)(\bR_i+\bB_j)(\bR_i+\frac{1}{2}(\bB_i+\bB_j))}\right]\\
	=&\bR_i\left[\frac{2\bR_i^2+2\bR_i(\bB_i+\bB_j)+\frac{1}{2}(\bB_i+\bB_j)^2-2\bR_i^2-2\bR_i(\bB_i+\bB_j)-2\bB_i\bB_j}{(\bR_i+\bB_i)(\bR_i+\bB_j)(\bR_i+\frac{1}{2}(\bB_i+\bB_j))}\right]\\
	=&\bR_i\left[\frac{\frac{1}{2}(\bB_i^2+\bB_j^2-2\bB_i\bB_j)}{(\bR_i+\bB_i)(\bR_i+\bB_j)(\bR_i+\frac{1}{2}(\bB_i+\bB_j))}\right]\\
	=&\bR_i\left[\frac{\frac{1}{2}(\bB_i-\bB_j)^2}{(\bR_i+\bB_i)(\bR_i+\bB_j)(\bR_i+\frac{1}{2}(\bB_i+\bB_j))}\right]
	\geq0.
\end{align*}

Now, consider a subset of nodes $V'\subset V$ such that $\bR_i=\bR_j$ for all $i,j\in V'$, and suppose $|V'|=n$. Without loss of generality, let $V'=\{1,\dots,n\}$. We assume the result holds for any subset of $n-1$ of these nodes. Thus, we can write that
\begin{align*}
	\sum_{i=1}^{n-1}P(Z_{i,1}=1) =\sum_{i=1}^{n-1}\frac{\bR_i}{\bR_i+\bB_i} \geq(n-1)\frac{\bR_1}{\bR_1+\frac{1}{n-1}\sum_{i=1}^{n-1}\bB_i}.
\end{align*}
Let $\bB=\frac{1}{n-1}\sum_{i=1}^{n-1}\bB_i$, and let $\bB^*_i=\frac{1}{n}((n-1)\bB+\bB_n)$, for all $i\in V'$. We note that $\bB^*_i=\frac{1}{n}\sum_{j\in V'}\bB_j$ for all $i\in V'$. We have
\begin{align}
	\sum_{i\in V'}P(Z_{i,1}=1)	&\geq(n-1)\frac{\bR_1}{\bR_1+\bB}+\frac{\bR_1}{\bR_1+\bB_n}\label{eq:eq7},\\
	\sum_{i\in V'}P^*(Z_{i,1}=1)	&=n\frac{\bR_1}{\bR_1+\frac{1}{n}((n-1)\bB+\bB_n)}\label{eq:eq8}.
\end{align}
To prove the result, we take the difference between the right-hand sides of~\eqref{eq:eq7} and~\eqref{eq:eq8}, obtaining
\begin{align*}
	&\bR_1\left[\frac{(n-1)(\bR_1+\bB_n)+(\bR_1+\bB)}{(\bR_1+\bB)(\bR_1+\bB_n)}-\frac{n}{\bR_1+\frac{1}{n}((n-1)\bB+\bB_n)}\right]\\
	=&\bR_1\left[\frac{n\bR_1+(n-1)\bB_n+\bB}{(\bR_1+\bB)(\bR_1+\bB_n)}-\frac{n}{\bR_1+\frac{1}{n}((n-1)\bB+\bB_n)}\right]\\
	=&\bR_1\left[\frac{(n\bR_1+(n-1)\bB_n+\bB)(\bR_1+\frac{1}{n}((n-1)\bB+\bB_n)-n(\bR_1+\bB)(\bR_1+\bB_n)}{(\bR_1+\bB)(\bR_1+\bB_n)(\bR_1+\frac{1}{n}((n-1)\bB+\bB_n))}\right]\\
	=&\bR_1\left[\frac{\frac{n-1}{n}(\bB^2+\bB_n^2-2\bB\bB_n)}{(\bR_1+\bB)(\bR_1+\bB_n)(\bR_1+\frac{1}{n}((n-1)\bB+\bB_n))}\right]\\
	=&\bR_1\left[\frac{\frac{n-1}{n}(\bB-\bB_n)^2}{(\bR_1+\bB)(\bR_1+\bB_n)(\bR_1+\frac{1}{n}((n-1)\bB+\bB_n))}\right]
	\geq0.
\end{align*}
\end{IEEEproof}

We now use this result to prove that for networks with symmetry, solutions to Problem~\ref{pb:init1} will be symmetric in nature as well. Suppose for a graph $\G=(V,\E)$ we have an automorphism $\sigma$ of $\G$. Assuming $\G$ is finite, we have that $Aut(\G)$ is a finite group, and thus we can use $\sigma$ to generate a cyclic subgroup $\langle\sigma\rangle$ of $Aut(\G)$. The subgroup $\langle\sigma\rangle$ is necessarily finite and of finite order. We also have that for any node $i$ in $V$ there exists some integer $k\in\integerspositive$ such that $i=\sigma^k(i)$. We will denote by $k_i$ the smallest such integer for which this holds for a given $i\in V$.




%
\begin{theorem}\label{lma:automorphism}
For a network $\G=(V,\E)$ with a given initialization $(\boldB,\boldR)$, consider an automorphism $\sigma$ of $\G$. Let $m$ be the order of the cyclic subgroup $\langle\sigma\rangle$ of $Aut(\G)$. Suppose $R_i=R_{\sigma(i)}$, for all $i\in V$. Let $\boldB^*$ be an alternate curing initialization such that $B^*_i=\frac{1}{m}\sum_{j=1}^{m}B_{\sigma^j(i)}$ for all $i\in V$. Then, $\tI_1^*\leq\tI_1$.
\end{theorem}

\begin{IEEEproof}
Note that for a given $i\in V$, $\bR_i=\sum_{j\in\N_i'}R_j$. We have that if $j\in \N_i'$ then $\sigma^k(j)\in\N_{\sigma^k(i)}'$ for any $k\in\integerspositive$. Since by assumption $R_j=R_{\sigma(j)}=R_{\sigma^k(j)}$ for all $j\in\N_i'$, we have that $\bR_i=\bR_{\sigma^k(i)}$ for all $i\in V$ and $k\in\integerspositive$. From Lemma~\ref{lma:lma2.5}, we need only to prove that $\bB^*_i=\frac{1}{k_i}\sum_{j=1}^{k_i}\bB_{\sigma^j(i)}$, for all $i\in V$. We start with
\begin{align*}
	\bB^*_i=\sum_{j\in\N_i'}B^*_j=\sum_{j\in\N_i'}\frac{1}{m}\sum_{l=1}^{m}B_{\sigma^l(j)}.
\end{align*}
We note that for a given positive integer $p$, 
\begin{align*}
	p\sum_{l=1}^{k_j}B_{\sigma^l(j)} = \sum_{l=1}^{k_j}(B_{\sigma^l(j)} + B_{\sigma^{l+k_j}(j)} + \dots + B_{\sigma^{l+(p-1)k_j}(j)}) = \sum_{l=1}^{pk_j}B_{\sigma^l(j)},
\end{align*}
which follows from the definition of $k_j$. Since $m$ is the order of $\langle\sigma\rangle$, we have $\sigma^m(j)=j$, and since $k_j$ is the smallest such integer for which this holds, we must have $k_j|m$ (i.e., $k_j$ divides $m$). Choosing $p=\frac{m}{k_j}$ and using the identity between the left-most and right-most terms in the above equation, we have that
\begin{align*}
	\frac{1}{k_j}\sum_{l=1}^{k_j}B_{\sigma^l(j)}=\frac{p}{pk_j}\sum_{l=1}^{k_j}B_{\sigma^l(j)}=\frac{1}{pk_j}\sum_{l=1}^{pk_j}B_{\sigma^l(j)}=\frac{1}{m}\sum_{l=1}^mB_{\sigma^l(j)},
\end{align*}
and thus
\begin{align*}
	\bB^*_i = &\sum_{j\in\N_i'}\frac{1}{m}\sum_{l=1}^mB_{\sigma^l(j)} = \frac{1}{m}\sum_{l=1}^m\sum_{j\in\N_i'}B_{\sigma^l(j)} = \frac{1}{m}\sum_{l=1}^m\sum_{j\in\N_{\sigma^l(i)}'}B_j \\
			= &\frac{1}{m}\sum_{l=1}^m\bB_{\sigma^l(i)} = \frac{1}{k_i}\sum_{l=1}^{k_i}\bB_{\sigma^l(i)}.
\end{align*}
The result then follows from Lemma~\ref{lma:lma2.5}.
\end{IEEEproof}

Another way to think about Theorem~\ref{lma:automorphism} is that the optimal allocation distributes resources equally within sets of the orbit partition by $\langle\sigma\rangle$. We assume that infection resources are distributed evenly within this orbit partition, a consequence of assuming that $R_i=R_{\sigma(i)}$ for all $i\in V$. We then prove redistributing curing resources within each set of the partition, i.e., letting $B^*_i=\frac{1}{m}\sum_{j=1}^{m}B_{\sigma^j(i)}$ for all $i\in V$, will result in a decrease in the overall average infection rate.

We finish this part with few remarks. First, due to the complicated nature of $\tI_n$, it is difficult to analytically prove this result for time $n>1$. Also, while Theorem~\ref{lma:automorphism} can be applied to a given automorphism, one can push this idea further. Consider a network $\G=(V,\E)$, and the orbit partition of $\G$ by $Aut(\G)$, denoted $P = \{V_1,\dots,V_q\}$, where $q\leq N$. For any given initialization, $(\boldB,\boldR)$, for which $\bR_i = \bR_j$  if $i$ and $j$ are similar under $Aut(\G)$, we postulate that we can reduce the average infection rate $\tI_n$ by redistributing curing resources evenly within each set of the partition. To elaborate, for each $i\in V$, if $i\in V_k$, we conjecture that taking $B_i^* = \frac{1}{|V_k|}\sum_{j\in V_k}B_j$ yields that $\tI_n(\boldB^*)\leq\tI_n(\boldB)$, for any $n$. This conjecture is supported empirically, but for reasons of space, we omit including a sample simulation.

\subsection{Gradient Descent}
As in~\cite{MH-FA-BG:18}, and further explored in Section~\ref{ch:curing-infection}, one can implement a gradient descent algorithm~\cite{DPB:99} to find solutions to Problem~\ref{pb:init1}. We first verify that this strategy can be employed for our initialization problem. For the curing problem presented in~\cite{MH-FA-BG:18}, the objective was to minimize the limiting average infection rate. As this measure is rather complicated to analyze, a proxy measure given by the expected average network exposure, $E[\tS_n|\F_{n-1}]$,  is used instead. For the initialization problem the equivalent method would be to study this value for time $n=1$, at which point this simply reduces to minimizing the average infection rate at time one, given by $\tI_1$, as a function of the curing initialization $\boldB$. In Problem~\ref{pb:init1}, we assign curing resources subject to a budget $\budgetB$. If for a general time $n$ we instead use $\tI_1$ as a simple proxy, it is easy to use gradient flow techniques to find an optimal policy for this one-step problem. We note that any optimal policy will make full use of the budget $\budgetB$.

\begin{proposition}\longthmtitle{Gradient descent conditions}\label{prop:grad_descent}
For a general network $\G=(V,\E)$ equipped with the Polya network contagion model, if for a given infection initialization $\boldR$ we have that $\bR_i>0$ for each $i\in V$, then the average infection rate at time one, $\tI_1$, is convex with respect to the curing initialization parameters $\boldB$. Moreover, the feasible set
\begin{align*}
	\X = \left\{\{B_i\}_{i=1}^N\in\realnonnegative^N \Bigg| \sum_{i=1}^NB_i=\budgetB\right\}
\end{align*}
is convex and compact.
\end{proposition}
The proof of this proposition is straightforward, with the matter of convexity of $\tI_1$ having been discussed earlier in this section. Of note, however, is the fact that we require $\bR_i>0$ for each $i\in V$, in order to ensure our feasible set $\X$ remains valid. If this is satisfied, we can then use a constrained gradient descent method~\cite[Chapter~2]{DPB:99}, as described in Algorithm~\ref{alg:grad_descent}. This algorithm will converge to an optimal solution for the one-step optimization problem, but not necessarily for the $n$-step case.
\begin{algorithm}[!t]
\caption{Constrained gradient descent on a simplex~\cite{DPB:99}}
\label{alg:grad_descent}
  \begin{algorithmic}
  	\STATE{$f(x_1,\dots,x_N)\gets$ function to be minimized}
	\\
	\STATE Start at an arbitrary node:
	\STATE $ y_1 = (y_{1,1}, \dots, y_{1,N}) = (\budgetB,0,\ldots,0) $
	\\
	\FOR{$ k = 1 : stoptime $}
		\STATE Find the direction of steepest descent:
		\STATE $ i = \argmin_{j \in V} \frac{\partial f}{\partial x_j} |_{y_k} $
		\\
		\STATE Move only in that direction:
		\STATE $ \bar{y}_{k,i} = \budgetB $, and $ \bar{y}_{k,j} = 0 $ for all $ j \neq i $
		\\
		\STATE Select the step size using the limit minimization rule:
		\STATE $ \alpha_k = \argmin_{\alpha \in [0,1]} f(y_k + \alpha(\bar{y}_k - y_k)) $
		\\
		\STATE Perform the gradient descent:
		\STATE $ y_{k+1} = y_k +  \alpha_k(\bar{y}_k - y_k) $
	\ENDFOR
  \end{algorithmic}
\end{algorithm}

\section{Curing-Infection Problems}\label{ch:curing-infection}

Another technique to control network infection for the Polya network model is to control the allocation of our curing parameters, $\{\Delta_{b,i}(n)\}_{i=1}^N$, subject to a fixed budget, $\budget$. As in the initialization setup, this can be done for a finite or infinite horizon. Previously,~\cite{MH-FA-BG:19} examined the problem of minimizing the limiting average infection rate and introduced a number of heuristic optimization policies. 
We begin by refocusing the problem in a game-theoretic setup.

We consider a two player game, where player one minimizes the average infection rate, $\tI_n$, while player two maximizes this same value. Player one controls the distribution of curing parameters, $\{\Delta_{b,i}(n)\}_{i=1}^N$, in accordance with a budget, $\budgetB$, while player two controls the distribution of the infection parameters, $\{\Delta_{r,i}(n)\}_{i=1}^N$  in accordance with a budget, $\budgetR$. The budget in each case is fixed for all time steps, and the allocation of resources for a given time $n$ is determined prior to any draws being made at that time. Thus, if either player allocates resources to a node for which the opposite colour ball is drawn, those resources will go to waste.

Formally, player one's objective is to find:
\begin{align}\label{eq:obj_curing}
	\min_{\{\{\Delta_{b,i}(k)\}_{i=1}^N:\sum_{i=1}^N\Delta_{b,i}(k)=\budgetB\},k=1,\ldots,n}\tI_n,
\end{align}
assuming the minimum exists, while player two's objective is to find:
\begin{align}\label{eq:obj_infect}
	\max_{\{\{\Delta_{r,i}(k)\}_{i=1}^N:\sum_{i=1}^N\Delta_{r,i}(k)=\budgetR\},k=1,\ldots,n}\tI_n,
\end{align}
assuming the maximum exists. We refer to \eqref{eq:obj_curing}-\eqref{eq:obj_infect} as the Delta-curing problem.

Since finding an optimal control policy for either player is not tractable for a general network,
we simplify the problem by adopting the expected network exposure, $E[\tS_n|\mathcal{F}_{n-1}]=E[\tS_n|\sigma(\{Z_j^{n-1}\}_{j=1}^N)]$, as the metric to optimize instead of the average infection rate. Using this (strongly correlated) proxy measure drastically facilitates the analysis.
We thus consider a two player zero-sum game, where player one minimizes the expected network exposure over the curing parameters $\{\Delta_{b,i}(n)\}_{i=1}^N$ and player two maximizes the same value over the infection parameters $\{\Delta_{r,i}(n)\}_{i=1}^N$. 
We begin with an important result about the expected network exposure.

\begin{proposition}\longthmtitle{Convexity-Concavity of Network Exposure}\label{thm:convex}
Let $\G = (V,\E)$ be a general network, and consider the Polya network contagion model on $\G$, with arbitrary initial conditions. Then, the expected network exposure $E[\tS_n|\F_{n-1}]$ is convex with respect to the curing parameters $\{\Delta_{b,i}(n)\}_{i=1}^N$ and concave with respect to the infection parameters $\{\Delta_{r,i}(n)\}^N_{i=1}$, for all $n$.
\end{proposition}
\begin{IEEEproof} 

Using~\eqref{eq:X_i,n} and~\eqref{eq:S_i,n}, we consider $E[\tS_n|\F_{n-1}]$ as a function of the vectors 
\begin{gather*}
	\bx=(x_1,\ldots,x_N)^T=(\Delta_{b,1}(n),\ldots, \Delta_{b,N}(n))^T\\
	\by=(y_1,\ldots,y_N)^T=(\Delta_{r,1}(n),\ldots, \Delta_{r,N}(n))^T, 
\end{gather*}
by reformulating~\eqref{eq:S_i,n} as follows:
\begin{align*}
	S_{i,n}=f_{i,n}(\bx,\by,Z_n) = \frac{c_i+\delta_i(\by,Z_n)}{c_i+d_i+\sigma_i(\bx,Z_n)+\delta_i(\by,Z_n)},
\end{align*}
where
\begin{align*}
	&c_i=\bR_i+\sum_{t=1}^{n-1}\sum_{j\in\N_i'}\Delta_{r,j}(t)Z_{j,t}, \quad d_i=\bB_i+\sum_{t=1}^{n-1}\sum_{j\in\N_i'}\Delta_{b,j}(t)(1-Z_{j,t}) \\
	&\delta_i(\by,Z_n)=\sum_{j\in\N_i'}y_j Z_{j,n}, \quad
	\sigma_i(\bx,Z_n)=\sum_{j\in\N_i'}x_j(1-Z_{j,n}).
\end{align*}

Alternatively, we let $A$ represent the adjacency matrix, including self-loops, for our given network $\G$, where $A_{ij}=1$ if and only if $i\in\N_j'$, and $A_{ij}=0$ otherwise. We then construct an $N\times N$ square matrix, $D$, where $D_{ij}=A_{ij}(1-Z_{j,n})$. Letting $D_i$ represent the $i^{th}$ row of the matrix $D$, we have that $\sigma_i(x,Z_n) = D_i\bx$.
Likewise, we can construct an $N\times N$ square matrix, $C$, where $C_{ij}=A_{ij}Z_{j,n}$. Letting $C_i$ represent the $i^{th}$ row of the matrix $C$, we have that $ \delta_i(y,Z_n) = C_i\by  $, and hence
\begin{align*}
	f_{i,n}(\bx,\by,Z_n) = \frac{c_i+C_i\by}{c_i+d_i+C_i\by+D_i\bx},
\end{align*}
where we have dropped the dependencies of the matrices $C$ and $D$ on $Z_n$ for simplicity.

Taking the expectation of $S_{i,n}$ given the history of the contagion process up to time $n-1$, we have
\begin{align}\label{eq:convex}
	E[S_{i,n}|\F_{n-1}]&=E[f_{i,n}(\bx,\by,Z_n)|\F_{n-1}] \nonumber\\
	&=\sum_{z_n\in\{0,1\}^N}f_{i,n}(\bx,\by,z_n)P(Z_n=z_n|Z^{n-1}=z^{n-1}).
\end{align}

We note that $P(Z_n=z_n|Z^{n-1}=z^{n-1})$ is independent of  our choice of $\bx$, and for any fixed realization $z_n = (z_{1,n},\ldots,z_{N,n})\in\{0,1\}^N$, $f_{i,n}(\bx,\by,z_n)$ is convex in $\bx$ over $\realnonnegative^N$ (the proof is given in~\cite{MH-FA-BG:19}). Hence, it follows that $E[S_{i,n}|\F_{n-1}]$ is convex in $\bx$ over $\realnonnegative^N$, and thus so too is $E[\tS_n|\F_{n-1}]$. Concavity in $\by$ follows from a symmetry argument since 
\begin{align*}
	1-f_{i,n}(\bx,\by,Z_n)= \frac{d_i+D_i\bx}{c_i+d_i+C_i\by+D_i\bx}
\end{align*}
is convex in $\by$, thus showing $f_{i,n}(\bx,\by,Z_n)$ is concave in $\by$. The rest of the proof follows as in the convexity argument for $\bx$.
\end{IEEEproof}

The natural symmetry of this model makes this result rather intuitive, and allows us to establish a nice setup for the game theoretic problem. We have proven that $E[\tS_n|\F_{n-1}]$ is convex in $\bx$ and concave in $\by$ over $\realnonnegative^N\times\realnonnegative^N$. In order to get a better understanding of this function, we consider the partial derivatives with respect to $\bx$ and $\by$. We note that
\begin{gather*}
	\nabla_{\bx}f_{i,n}(\bx,\by,z_n)=\frac{-\nabla_{\bx}D_i\bx(c_i+C_i\by)}{(c_i+d_i+C_i\by+D_i\bx)^2}\\
	\nabla_{\by}f_{i,n}(\bx,\by,z_n)=\frac{\nabla_{\by}C_i\by(d_i+D_i\bx)}{(c_i+d_i+C_i\by+D_i\bx)^2},
\end{gather*}
and furthermore,
\begin{align*}
	\frac{\partial}{\partial x_j}D_i\bx=A_{ij}(1-z_{j,n}),\ \mathrm{and} \
	\frac{\partial}{\partial y_j}C_i\by=A_{ij}z_{j,n}.
\end{align*}
Since $A_{ij}$ and $z_{j,n}$ can only take values in $\{0,1\}$, we have that $\frac{\partial}{\partial x_j}f_{i,n}(\bx,\by,z_n)\leq0$ and $\frac{\partial}{\partial y_j}f_{i,n}(\bx,\by,z_n)\geq0$. Using this fact in conjunction with~\eqref{eq:convex}, we determine that over the space $\X\times\Y=\mathbb{R}^N_{\geq0}\times\mathbb{R}^N_{\geq0}$, $E[\tS_n|\F_{n-1}]$ has no saddle point. 
Given our fixed allocation budgets $\budgetB$ and $\budgetR$, we restrict ourselves to considering sets of the form $\X\times\Y = \{\{\Delta_{b,i}(n)\}_{i=1}^N\in\realnonnegative^N|\sum_{i=1}^N\Delta_{b,i}(n)\leq\budgetB\}\times\{\{\Delta_{r,i}(n)\}_{i=1}^N\in\realnonnegative^N|\sum_{i=1}^N\Delta_{r,i}(n)\leq\budgetR\}$.
Returning to our game over the expected network exposure, we remark that for any given $n$, the sets $\X = \{\{\Delta_{b,i}(n)\}_{i=1}^N\in\realnonnegative^N|\sum_{i=1}^N\Delta_{b,i}(n)\leq\budgetB\}$ and $\Y = \{\{\Delta_{r,i}(n)\}_{i=1}^N\in\realnonnegative^N|\sum_{i=1}^N\Delta_{r,i}(n)\leq\budgetR\}$ are convex and compact. This gives rise to the following result.

\begin{theorem}\longthmtitle{Nash Equilibrium for Network Exposure}\label{thm:nash_exposure}
Let $\G=(V,\E)$ be a general network equipped with the Polya network contagion model, with arbitrary initial conditions. For a given time $n$, consider a two-player zero-sum game where player one tries to minimize $E[\tS_n|\F_{n-1}]$ over the parameters $\{\Delta_{b,i}(n)\}_{i=1}^N$ and player two tries to maximize $E[\tS_n|\F_{n-1}]$ over the parameters $\{\Delta_{r,i}(n)\}_{i=1}^N$. Then, if we take our set of allowable policies to be of the form $\X\times\Y = \{\{\Delta_{b,i}(n)\}_{i=1}^N\in\realnonnegative^N|\sum_{i=1}^N\Delta_{b,i}(n)\leq\budgetB\}\times\{\{\Delta_{r,i}(n)\}_{i=1}^N\in\realnonnegative^N|\sum_{i=1}^N\Delta_{r,i}(n)\leq\budgetR\}$, the resulting game admits a Nash equilibrium . Moreover, the equilibrium policy will satisfy $\sum_{i=1}^N\Delta_{b,i}(n)=\budgetB$ and $\sum_{i=1}^N\Delta_{r,i}(n)=\budgetR$.
\end{theorem}
\begin{IEEEproof}
Since the function is convex-concave and over a compact set, the existence follows from the classical minimax theorem, see~\cite{DPB-AN-AEO:03}. By the definition of $E[\tS_n|\F_{n-1}]$, and since the function has no saddle point in the interior of its domain, the optimal policy will utilize the full budget.
\end{IEEEproof}

The equilibrium policy from Theorem~\ref{thm:nash_exposure} can be determined numerically using gradient descent algorithms~\cite{DPB:99}, as we are optimizing a convex/concave function over a simplex, but for large networks such algorithms can be computationally expensive, especially when considering the complexity of $E[\tS_n|\F_{n-1}]$ as seen in~\eqref{eq:convex}. 
We next develop efficient heuristic strategies for both the initialization and Delta-curing problems
(studied in Section~\ref{ch:initialization} and the current section, respectively). Naturally, their effectiveness
will be measured via the original metric intended for the problems, which is average network infection rate.
\section{Heuristic Strategies}\label{sec:heuristics}
We provide here heuristic strategies based on the analysis performed thus far, while also integrating some of the previously tested techniques in~\cite{MH:17},~\cite{MH-FA-BG:19}.

When developing heuristics in the initialization case, we make heavy use of Theorem~\ref{lma:316}, which proved that any optimal solution to Problem~\ref{pb:init1} would allocate no resources to outer nodes (i.e., nodes with a nested neighbourhood). For the general Delta-curing problem, we have no such result, as it is possible to manufacture situations where the equilibrium policy as in Theorem~\ref{thm:nash_exposure} will assign resources to outer nodes. Such cases arise under specific circumstances, and it is uncommon for these nodes to receive significant resources. Thus, we introduce a similar set of resource allocation policies for both the initialization and Delta-curing problems of Sections~\ref{ch:initialization} and~\ref{ch:curing-infection}, respectively. While any of these resource allocation policies can be similarly applied to a player controlling the distribution of infection resources, for the purposes of our analysis here, we will consider these strategies from the perspective of the player controlling the curing resources.

\subsection{Interior Node Targeting}\label{subsec:heur_inner}
We now discuss our first strategy. Given a network, we determine which nodes are inner nodes, i.e., which nodes have a neighbourhood that is not a strict subset of another node's neighbourhood. We then distribute initialization resources uniformly among these nodes. If we let $\budgetB$ represent our curing initialization budget, and let $V'=\{i\in V|\N_i'\not\subset\N_j',\forall j\in V\}$, we allocate our resources according to the following scheme:
\begin{align}\label{eq:unif-balls}
	B_i=\begin{cases}
		\frac{1}{|V'|}\budgetB	&\text{ if $i\in V'$}\\
		0						&\text{ else}.
	\end{cases}
\end{align}
For the Delta-curing problem, a similar policy is employed as follows:
\begin{align}\label{eq:unif-curing}
	\Delta_{b,i}(n)=\begin{cases}
		\frac{1}{|V'|}\budgetB	&\text{ if $i\in V'$}\\
		0						&\text{ else},
	\end{cases}
\end{align}

We can also combine these schemes with the heuristic strategies in~\cite{MH:17},~\cite{MH-FA-BG:19}. Among these strategies, the best performing heuristic strategy allocates resources based on a combination of centrality measures,  node degree and closeness centrality, as well as the infection level of each node's super urn. In the case of initialization policies, as they are implemented before the onset of contagion, we implement an allocation strategy based on only the node degree and closeness centrality.

The degree of a node is simply measured as the size of a node's neighbourhood, $|\N_i|$. The closeness centrality of a node, $C_i$, is a measure used to determine a node's topological position within a network, defined as
\begin{align}\label{eq:closeness}
	C_i=\frac{1}{\sum_{j\in V}d(i,j)},
\end{align}
where $d(i,j)$ is the shortest path length from node $i$ to node $j$. Thus the node with the highest closeness centrality has the shortest average path length to every node within the network. We can distribute resources as follows:
\begin{align}\label{eq:centrality-balls}
	B_i=\begin{cases}
		\frac{|\N_i|C_i}{\sum_{j\in V'}|\N_j|C_j}\budgetB	&\text{ if $i\in V'$}\\
		0												&\text{ else}.
	\end{cases}
\end{align}

In the Delta-curing setup, we also account for the evolution of the super urn proportion, $S_{i,n-1}$, of each node. As before, we let $C_i$ represent the closeness centrality of a given node $i$, calculated by~\eqref{eq:closeness}. We then allocate resources as follows:
\begin{align}\label{eq:centrality-curing}
	\Delta_{b,i}(n)=\begin{cases}
		\frac{|\N_i|C_iS_{i,n-1}}{\sum_{j\in V'}|\N_j|C_jS_{j,n-1}}\budgetB	&\text{ if $i\in V'$}\\
		0																&\text{ else}.
	\end{cases}
\end{align}
Alternatively, similar strategies based on ``betweenness centrality" can be used~\cite{MH-FA-BG:19}.

\subsection{Minimized Node Targeting}\label{subsec:heur_fewer}
We can improve upon our initialization and curing policies by using the idea of minimizing the number of nodes we target at once. We introduce a method that attempts to target nodes in such a manner that all nodes have resources within their immediate neighbourhood $\N_i'$ (i.e., ensuring $\bB_i>0$, $\forall i\in V$), while targeting the fewest possible nodes. We use Algorithm~\ref{alg:node_target} to determine which nodes will be allocated resources.

This algorithm starts by locating outer nodes in the network and then targets nodes directly adjacent these outer nodes. All nodes that are within the neighbourhood of the targeted nodes are then removed from the set of nodes being evaluated. Within this new, smaller set of nodes, we once again locate outer nodes and repeat the process. In the case that there are no outer nodes within the set of nodes being evaluated, all remaining nodes are added to the set of targeted nodes.
The purpose of this algorithm is to effectively employ Theorem~\ref{lma:316} and ensure that our curing policy does not target outer nodes, while still requiring that $\bB_i>0$ for these nodes. Once the set of outer nodes is excluded from directly receiving resources, the algorithm then works inwards, attempting to minimize the number of targeted nodes.

\begin{algorithm}[!t]
\caption{Node Targeting Algorithm}
\label{alg:node_target}
	\begin{algorithmic}
		\STATE $V \gets$ set of all nodes
		\STATE $V' \gets$ target set
		\STATE $V_{test} \gets$ nodes to test for this loop
		\\
		\STATE Start with $V'=\emptyset,$ $V_{test}=V$
			\WHILE{$V_{test}\neq\emptyset$}
				\STATE Identify outer nodes:
				\STATE Set $V_{outer}=\{i\in V_{test}|\N_i'\subset\N_j',\text{ for some } j\in V_{test}\}$
					\IF{$V_{outer}=\emptyset$}
						\STATE $V' = V'\cup V_{test}$
						\STATE \textbf{break}
					\ENDIF
				\STATE Target nodes adjacent to outer nodes:
				\STATE Set $V_{inner}=V_{test}\backslash V_{outer}$
				\STATE Set $V_{added\_targets}=\{i\in V_{inner}|i\in\N_j'\text{ for some }j\in V_{outer}\}$
				\STATE Set $V'=V'\cup V_{added\_targets}$
				\STATE Update parameters for next loop:
				\STATE Set $V_{coverage}=\{i\in V_{test}|i\in\N_j'\text{ for some }j\in V_{added\_targets}\}$
				\STATE Set $V_{test}=V_{test}\backslash V_{coverage}$
			\ENDWHILE
	\end{algorithmic}
\end{algorithm}
Once the set of nodes to be targeted $V'$ has been determined, we can once again distribute resources in one of two manners as in the previous strategy set. We can distribute resources to these nodes uniformly according to~\eqref{eq:unif-balls}.
Alternatively, we can distribute resources to more central nodes among those targeted
using~\eqref{eq:centrality-balls}.
We note that these heuristics are entirely based on network structure, and by design will symmetrically allocate resources in accordance with Theorem~\ref{lma:automorphism}. Furthermore, by limiting the number of nodes that are targeted, we limit the computational complexity of determining resource allocation. 

We can apply the same techniques for the Delta-curing setup, where we consider a uniform allocation within the target set via~\eqref{eq:unif-curing},
or use a strategy that accounts for the centrality and super urn proportion of the targeted nodes according to~\eqref{eq:centrality-curing}.

\subsection{Dense Networks}\label{subsec:heur_dense}


We now consider a limitation of the heuristics presented thus far. Our previous sets of solutions depend on the assumption that a network will possess a set of nodes whose respective neighbourhoods are strict subsets of some other node's neighbourhood. 
While this may be a useful tactic for sparse networks, it is a stringent requirement for dense networks. 
To account for this, we introduce a final node targeting strategy that does not have this limitation, but works based on similar principles.

Our previous methods 
utilized the idea of leveraging Theorem~\ref{lma:316}, in order to prevent suboptimal node targeting while still ensuring $\bB_i>0$, for every node $i \in V$. We introduce an alternative method for defining inner and outer nodes when the conditions of Theorem~\ref{lma:316} are not met.

This technique, presented in Algorithm~\ref{alg:node_target_dense}, begins by ranking the nodes by a centrality measure (closeness centrality for our purposes), from highest to lowest. Nodes are added one-by-one to the set of targeted nodes, $V'$, in order of rank. Once the set of targeted nodes provides total coverage of the network (i.e., the neighbourhoods of these nodes contain the entire network), the algorithm stops. For an additional reduction, we go through the list of targeted nodes in reverse order and remove any nodes not required to maintain full coverage of the network. 

\begin{algorithm}[!t]
\caption{Node Targeting Algorithm for Dense Networks}
\label{alg:node_target_dense}
	\begin{algorithmic}
		\STATE $V \gets$ set of all nodes
		\STATE $V' \gets$ nodes to be targeted
		\STATE $C_i \gets$ centrality score of given node
		\\
		\STATE Without loss of generality, assume $C_1>C_2>\dots>C_N$, where $N$ is the network size
		\\
		\STATE Start with $V'=\emptyset,$ $V_{coverage}=\emptyset$, n=0
			\WHILE{$V_{coverage}\neq V$}
				\STATE Add most central node to target set:
				\STATE $V'=V'\cup \{n\}$
				\STATE Update loop parameters:
				\STATE Set $V_{coverage}=\{i\in V|i\in\N_j',\text{ for some }j\in V'\}$
				\STATE n++
			\ENDWHILE
		\\
		\STATE \textit{Optional}:
			\FOR{$n = |V'| : 1$}
				\STATE Let $V_{coverage} = \{i\in V|i\in\N_j',\text{ for some }j \in V'\backslash\{n\}\}$
				\IF{$V_{coverage}=V$}
					\STATE Set $V' = V'\backslash\{n\}$
				\ENDIF
			\ENDFOR
	\end{algorithmic}
\end{algorithm}

Once the target set has been determined, we can distribute resources either uniformly among these nodes via~\eqref{eq:unif-balls},
or preferentially according to centrality as per~\eqref{eq:centrality-balls}.
Likewise for the Delta-curing case, we can allocate resources uniformly among these nodes
using~\eqref{eq:unif-curing}
or based on centrality and super urn proportion of the targeted nodes via~\eqref{eq:centrality-curing}.

\section{Simulation Results}\label{ch:simulations}

We evaluate the performance of the heuristic policies detailed in Section~\ref{sec:heuristics}, and compare their performance to supposed `optimal' strategies obtained through gradient descent (note that we say `optimal' because the gradient descent only converges to an optimal policy for the one-step proxy measures we introduced in both the initialization and Delta-curing cases). Algorithm~\ref{alg:polya_process} provides the general format of the simulations, while specifics of both the initialization and Delta-curing setups are provided in Sections~\ref{sec:init_trials} and~\ref{sec:curing_trials}, respectively.

\begin{algorithm}[!t]
\caption{Simulation Setup}
\label{alg:polya_process}
	\begin{algorithmic}
	  	\STATE $ A \gets $ adjacency matrix of the network
		\STATE $ \textit{numTrials} \gets $ trials to run for given policy
		\STATE $ \textit{steps} \gets $ number of time steps for each trial
		\STATE $ \textit{initPolicy} \gets $ initialization policy (uniform in no control case)
		\STATE $ \textit{curePolicy} \gets $ Delta-curing policy (uniform in no control case)
		\\
		\STATE Assign $ (\boldB,\boldR) $ using \textit{initPolicy} (under initialization budget)
		\FOR{s = 1 : \textit{numTrials}}
			\STATE $\vec{Z}_{s} \gets $ \textsc{RunTrial}($A, \boldR, \boldB, \budgetR, \budgetB, \textit{steps}, \textit{curePolicy}$)
		\ENDFOR
	  	\\
		\STATE{RunTrial}{$A, \boldR, \boldB, \budgetR, \budgetB, \textit{steps}, \textit{curePolicy}$}
			\STATE Initialize $ S_{i,0}$ using $ R_i $ and $ B_i $ for all $ i \in V $
			\FOR{t = 1 : \textit{steps}}
				\STATE Assign $ \Delta_{b,i}(t) $, $ \Delta_{r,i}(t) $ using \textit{curePolicy} (under curing budget)
				\STATE Generate $ \vec{Y} \sim \mathbf{Uniform}([0,1])^{N} $
				\IF{$ Y_i \leq S_{i,t-1} $}
					\STATE $ Z_{i,t} = 1 $
				\ELSE
					\STATE $ Z_{i,t} = 0 $
				\ENDIF
				\STATE Update $ S_{i,t} $ using $ \Delta_{r,i} $ and $ \Delta_{b,i}(t) $ for all $ i \in V $
			\ENDFOR
	\end{algorithmic}
\end{algorithm}

\begin{figure}[!ht]
\centering{
  	\subfigure[Barabasi-Albert network with 100 nodes and 99 edges.]{ \includegraphics[width=0.4\linewidth]{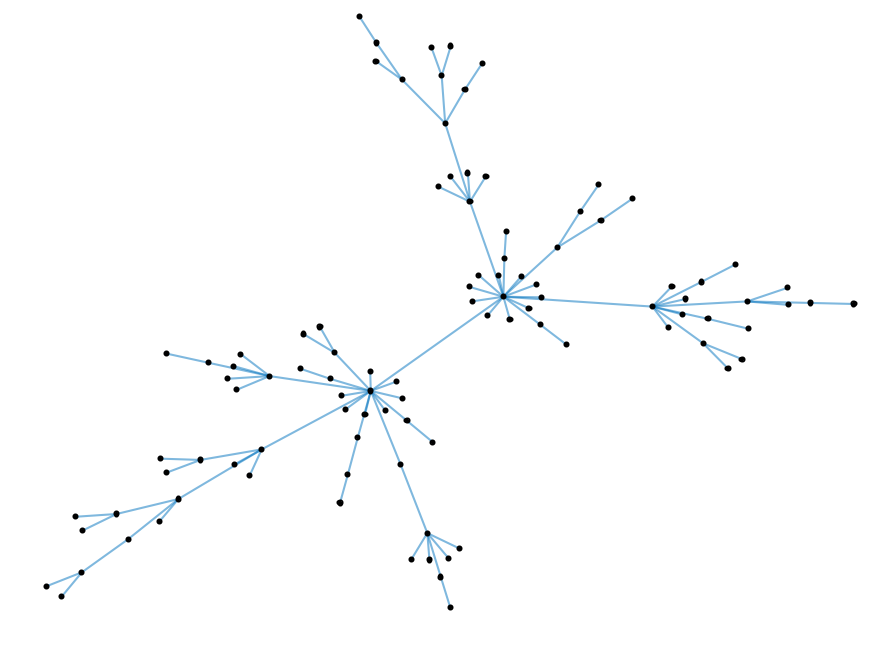}\label{fig:BANet}}
}
{
  	\subfigure[Barabasi-Albert network with 100 nodes and 945 edges.]{ \includegraphics[width=0.4\linewidth]{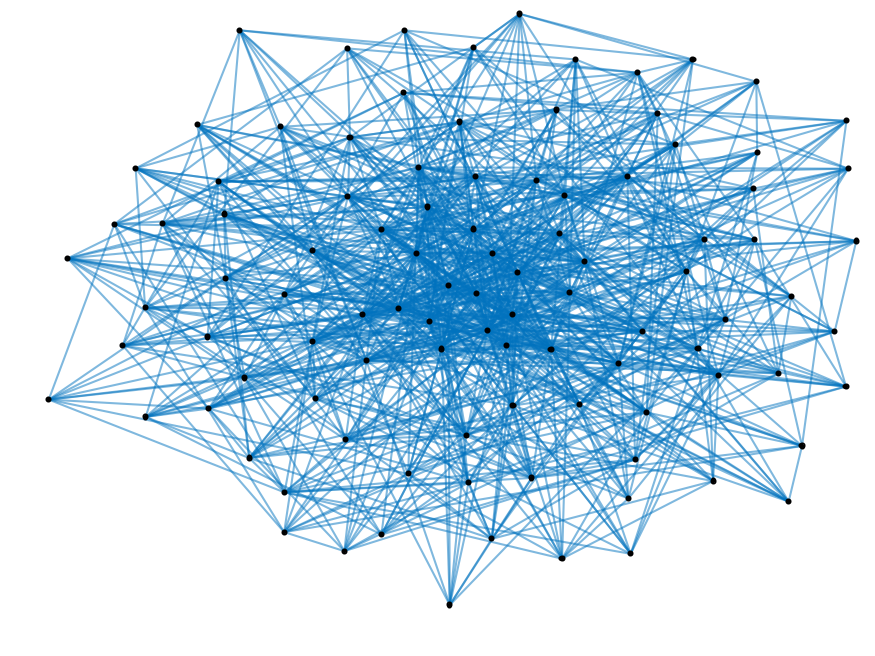}\label{fig:BADNet}}
}
\\
{
  	\subfigure[Facebook network with 1,363 nodes and 2,425 edges.]{ \includegraphics[width=0.5\linewidth]{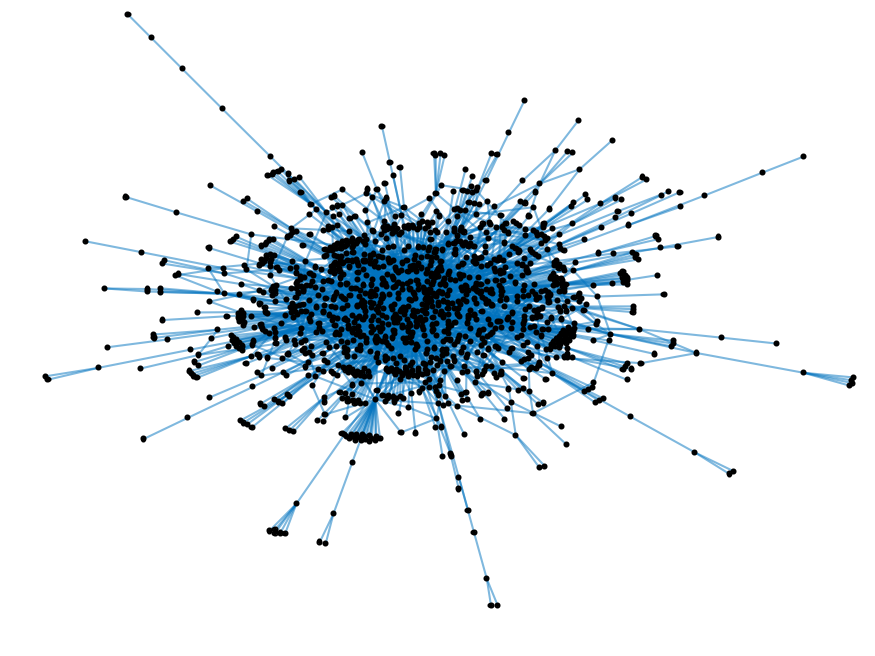}\label{fig:FaNet}}
}
\caption[Networks used for simulations]{Mid-to-large scale networks used for simulation purposes. Adjacency matrices for each of these networks is available online at \href{https://bit.ly/2ygLEqg}{\color{blue}{https://bit.ly/2ygLEqg}}.}
\label{fig:Networks}
\end{figure}

In order to ensure effectiveness of different strategies, we test these policies for a few different networks, each of which is depicted in Figure~\ref{fig:Networks}. These networks include algorithmically generated 100 node Barabasi-Albert networks, depicted in Figures~\ref{fig:BANet} and~\ref{fig:BADNet}, with average densities of 1.98 and 18.9 connections per node, respectively. The third network, depicted in Figure~\ref{fig:FaNet}, was generated by a tool~\cite{BR:13} that crawled through posts in a Facebook group in order to establish connections between users based on interactions on user posts (either via commenting, or `likes' on the post or comments). The resulting 1,363 node graph contains an average density of 3.56 connections per node, and provides a topology for a real-world social network.

\subsection{Targeting Algorithms}

Since the majority of the following analysis applies curing policies based on targeting subsets of nodes, we provide a brief comparison of the different node targeting algorithms. 

\begin{figure}[!ht]
\centering{
  	\subfigure[Example set of inner nodes.]{ \includegraphics[width=0.4\linewidth]{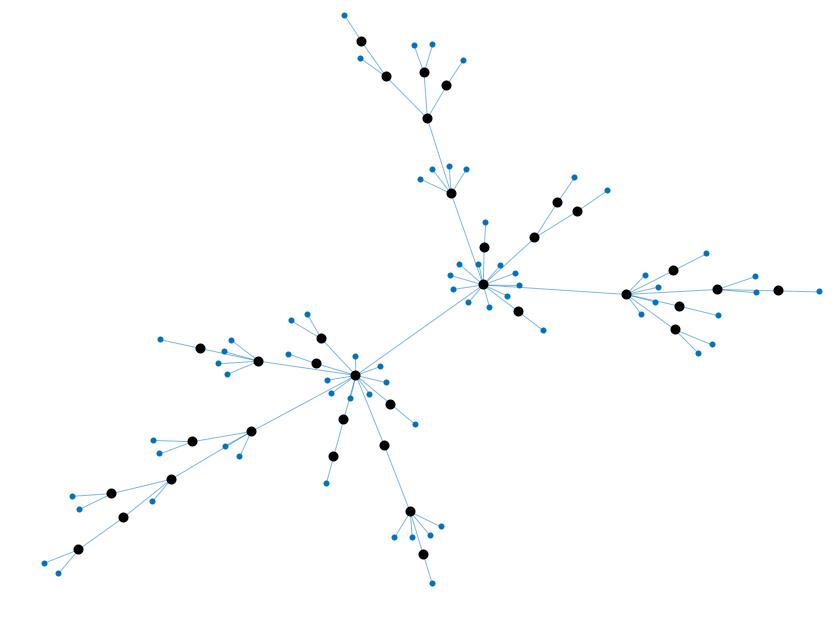}\label{fig:BANetInner}}
}
{
  	\subfigure[Example target set using Algorithm~\ref{alg:node_target}.]{ \includegraphics[width=0.4\linewidth]{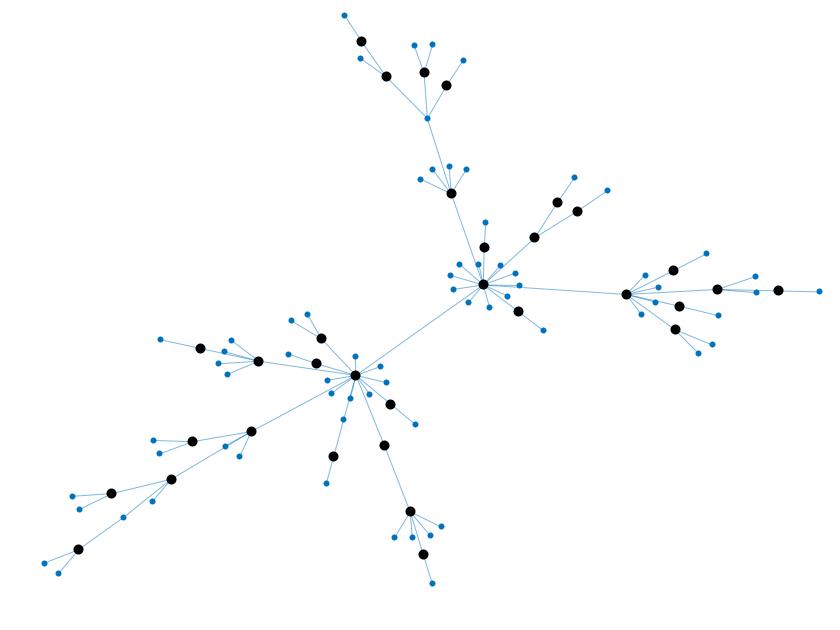}\label{fig:BANetAlg2}}
}
\\
{
  	\subfigure[Example target set using Algorithm~\ref{alg:node_target_dense} for sparse networks.]{ \includegraphics[width=0.4\linewidth]{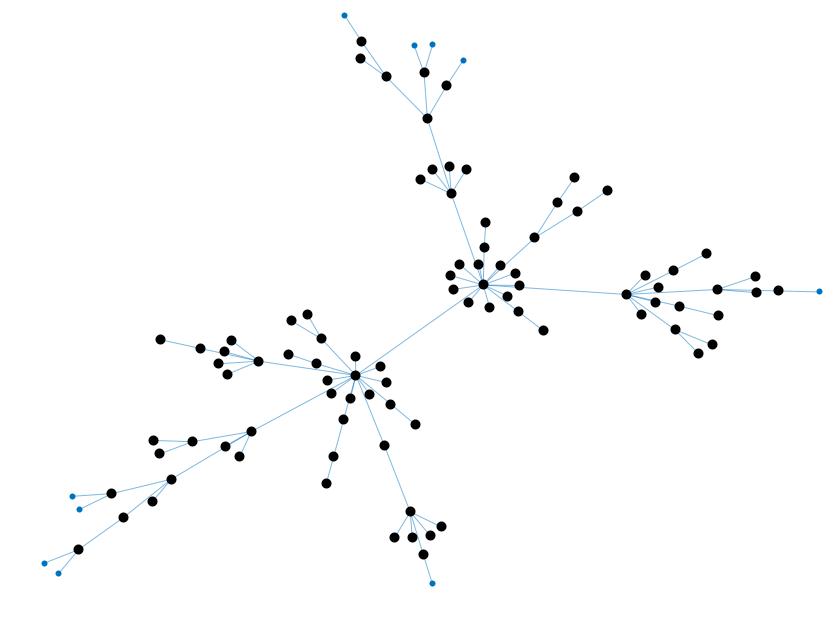}\label{fig:BANetAlg3}}
}
{
  	\subfigure[Example target set using Algorithm~\ref{alg:node_target_dense} for dense networks.]{ \includegraphics[width=0.4\linewidth]{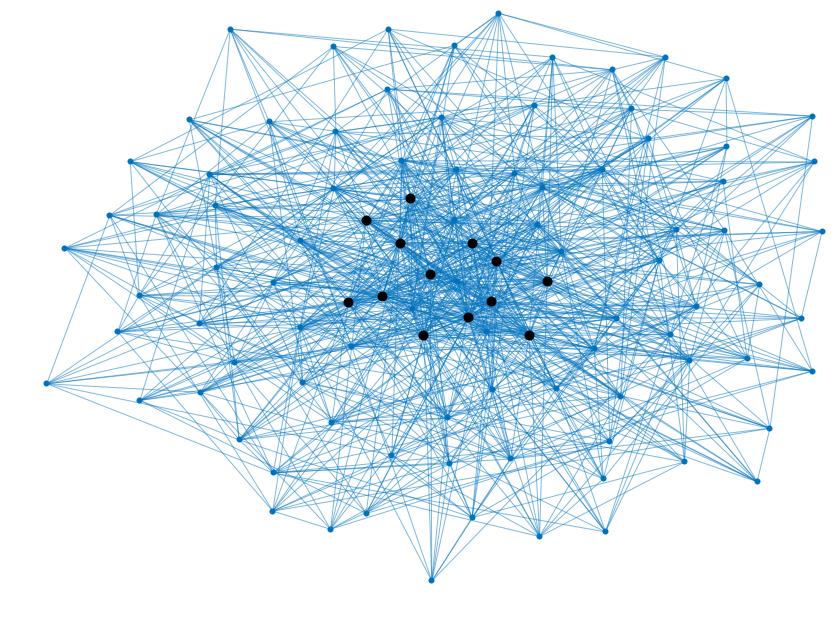}\label{fig:BADNetAlg3}}
}
\caption[Comparison of node targeting algorithms]{Comparison of target sets for various node targeting algorithms. In each case, the set of targeted nodes is enlarged and coloured black. It is important to note that Figure~\ref{fig:BANetAlg3}, as depicted here, does not apply the optional component of Algorithm~\ref{alg:node_target_dense}. In fact, if we instead choose to apply that option, we end up getting the same target set as shown in Figure~\ref{fig:BANetAlg2}.}
\label{fig:NodeTargets}
\end{figure}

Figure~\ref{fig:NodeTargets} depicts a comparison of the different targeting algorithms for the sparse network seen in Figure~\ref{fig:BANet}. For this specific network, the target set in Figure~\ref{fig:BANetAlg2} is a subset of the set of inner nodes in Figure~\ref{fig:BANetInner}, which itself is a subset of the target set shown in Figure~\ref{fig:BANetAlg3}. While, this is not a general rule, it is a common trend for very low-density networks. Algorithm~\ref{alg:node_target_dense} is not designed for use over sparse networks, and does not successfully leverage Theorem~\ref{lma:316}. Thus, we know that distributing resources among this set of target nodes will result in inefficiencies.
If we compare this to how Algorithm~\ref{alg:node_target_dense} applies to denser networks, as shown in Figure~\ref{fig:BADNetAlg3}, we observe a significant reduction in size of the set of target nodes. Algorithm~\ref{alg:node_target}, on the other hand, will be entirely ineffective at determing a subset of nodes to target for this same network, due to the lack of nested nodes.

\subsection{Initialization Trials}\label{sec:init_trials}

We start by comparing solutions to Problem~\ref{pb:init1}, and evaluating their effectiveness for various network compositions. Table~\ref{tab:init_strats} provides an overview of the considered initialization strategies.
More detailed explanations of these strategies can be found in Sections~\ref{ch:initialization} and~\ref{sec:heuristics}. It is important to note that we do not include any trials that use the optional component of Algorithm~\ref{alg:node_target_dense}, simply because we have the same target set as Algorithm~\ref{alg:node_target} when applied to the networks in~\ref{fig:BANet} and~\ref{fig:FaNet}, and it does not provide any further reduction to the target set when applied to the network in Figure~\ref{fig:BADNet}. The same reasoning applies to our curing strategies in Table~\ref{tab:curing_strats}.

\begin{table}[t!]
	\centering
	\caption{Initialization Strategies}\label{tab:init_strats}
	\begin{tabular}{cc}\hline\hline	
		$ \text{(i)} $ & Constrained gradient descent on a simplex: \\
			& Find $B_i $ using Algorithm~\ref{alg:grad_descent} on the function $\tI_1$ \\
\hline
		$ \text{(ii)} $ & Uniformly allocate the budget to all nodes: \\
			& $ B_i = \frac{\budgetB}{N} $ \\
	\hline
		$ \text{(iii)} $ & Uniformly allocate the budget to inner nodes: \\
			& $ B_i=\begin{cases}
						\frac{1}{|V'|}\budgetB	&\text{ if $i\in V'$} \\
						0						&\text{ else} \\
			\end{cases} $, where $V' = \{k\in V|\N_k'\not\subset\N_j',\forall j\in V\}$ \\
	\hline
		$ \text{(iv)} $ & Ratio of degree and closeness centrality  (inner): \\
			& $	B_i=\begin{cases}
					\frac{|\N_i|C_i}{\sum_{j\in V'}|\N_j|C_j}\budgetB	&\text{ if $i\in V'$} \\
					0												&\text{ else}
	\end{cases} $ \\
	\hline
		$ \text{(v)} $ & Allocate as in $ \text{(iii)} $ using $V'$ from Algorithm~\ref{alg:node_target} \\
	\hline
		$ \text{(vi)} $ & Allocate as in $ \text{(iv)} $, using $V'$ from Algorithm~\ref{alg:node_target} \\
	\hline
		$ \text{(vii)} $ & Allocate as in $ \text{(iii)} $ using $V'$ from Algorithm~\ref{alg:node_target_dense} (without Optional) \\
	\hline
		$ \text{(viii)} $ & Allocate as in $ \text{(iv)} $, using $V'$ from Algorithm~\ref{alg:node_target_dense} (without Optional) \\
	\hline
		$ \text{(ix)} $ & Ratio of degree and closeness centrality: \\
			& $ B_i=\frac{|\N_i|C_i}{\sum_{j\in V}|\N_j|C_j}\budgetB $ \\
\hline\hline
	\end{tabular}
\end{table}

\begin{figure}[!hbt]
\centering{
  	\subfigure[Results of initialization trials for low density Barabasi-Albert network depicted in~\ref{fig:BANet}]{\includegraphics[width=0.7\linewidth]{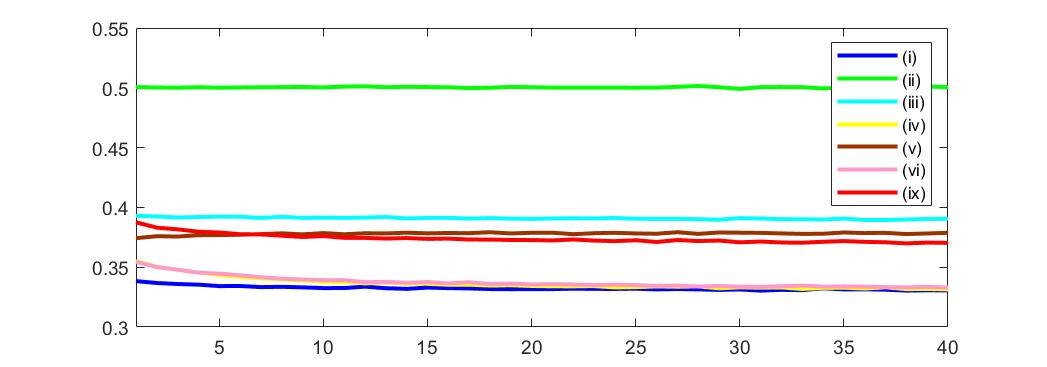}\label{fig:BAInitTrials}}
}
\\
{
  	\subfigure[Results of initialization trials for high density Barabasi-Albert network depicted in~\ref{fig:BADNet}]{\includegraphics[width=0.7\linewidth]{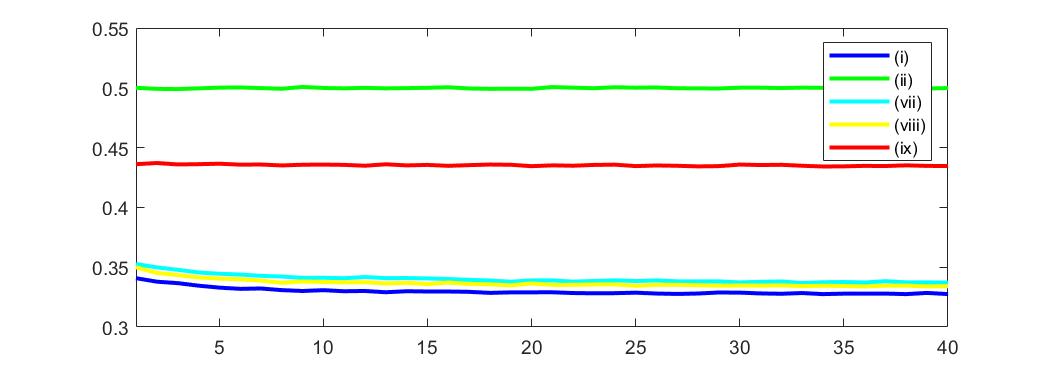}\label{fig:BADInitTrials}}
}
\\
{
  	\subfigure[Results of initialization trials for Facebook network depicted in~\ref{fig:FaNet}]{
  	\includegraphics[width=0.7\linewidth]{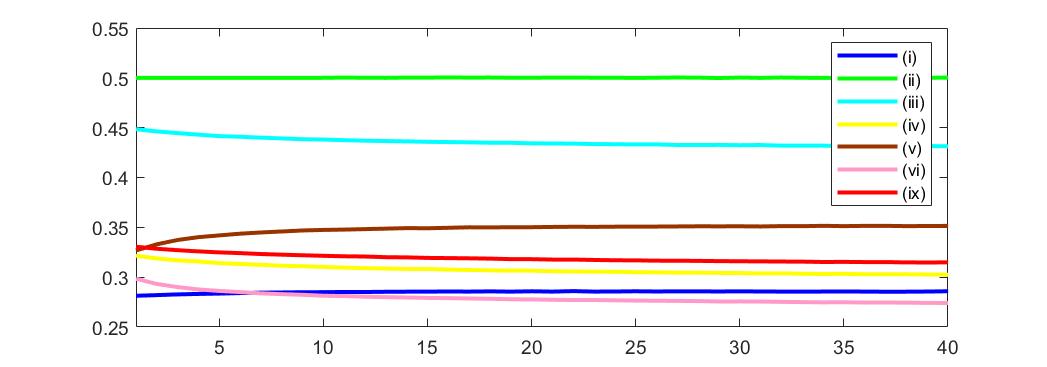}\label{fig:FaInitTrials}}
}
\caption[Comparison of initialization policies]{Plots of empirical average infection rate, $\tI_n$, for various initialization strategies across multiple networks. Lower values indicate lower levels of infection. In each case the network was given a fixed initialization budget of $\budgetR = \budgetB = 10N$, with $R_i=10$ $\forall i\in V$, and the Polya network contagion process was run over 1000 trials with a fixed $\Delta=5$.}
\label{fig:InitTrials}
\end{figure}

We follow the same process for each of the three networks in Figure~\ref{fig:Networks}. We begin with equal initialization budgets for both red and black resources (i.e., take $\budgetB=\budgetR$), and let the infection initialization be uniform over the network (i.e., take $R_i=\frac{\budgetR}{N}$, for all $i \in V$). We choose our curing initialization $\boldB$ in accordance with the desired strategy from Table~\ref{tab:init_strats}. Then, we run the Polya network contagion process with $\Delta_{r,i}(n)=\Delta_{b,i}(n)=\Delta$ fixed for each node $i\in V$ and each time $n$. We record the draw value, $Z_{i,n}$, for each node $i$ at each time $n$. The draw values are averaged over the number of trials that are run and then averaged over the set of nodes in the network to obtain an empirical measure for the average infection rate $\tI_n$ at a given time $n$. The resulting output of these trials can be seen in Figure~\ref{fig:InitTrials}.

We start our comparison with a few observations. First, as expected, the solution determined using the gradient descent algorithm results in the lowest empirical average infection rate at time $n=1$. For time $n>1$, however, this is not always the case. Figure~\ref{fig:FaInitTrials} is of particular interest, as we clearly see that the gradient descent strategy begins with a much lower infection rate, but over time is surpassed by strategy $(vi)$, which targets a reduced subset of nodes in accordance with Algorithm~\ref{alg:node_target}, with resources being focused more around central nodes.

The fact that the best strategy at time $n=1$ is not optimal as $n$ increases supports our hypothesis that the optimal solution will change for increasing values of $n$ and as $n\to\infty$. It is important to note, however, that our gradient algorithm merely converges to the optimal solution at time $n=1$, and thus the solution we have from using this algorithm is not necessarily optimal, as the algorithm has not been run indefinitely. Indeed, we note that had we chosen to either use a finer filter when implementing the limit minimization rule, as described in Algorithm~\ref{alg:grad_descent}, or increased the stopping time, or even chosen a step size that allowed for faster convergence, we would have seen an improved performance using this strategy (at least at time $n=1$).

Regardless, we consistently see that targeting a reduced subset of nodes, whether via Algorithms~\ref{alg:node_target} or~\ref{alg:node_target_dense}, results in performance that is close to or better than that obtained using gradient descent. These strategies, $(vi)$ and $(viii)$, always provide the best results of all the heuristic strategies being compared. This illustrates that we obtain near optimal performance while only distributing resources amongst a greatly reduced subset of nodes. Performance between the other heuristic strategies varies greatly depending on the network, though we note that even just targeting the inner nodes using strategy $(iii)$ sees a significant improvement over the no control case (i.e., the results obtained using a uniform distribution of initialization resources). 

We also remark that the general curve for each of these trials remains relatively flat. This is mostly due to the fact that $\Delta$ is assigned a fixed value. At a given time for a given node the expected number of red balls added to that particular node's super urn will be proportional to the average super urn proportion of the neighbouring nodes. This causes the super urn proportion of a given node to trend towards the average super urn proportion of its neighbours, which often takes a similar value because these nodes have overlapping sets of neighbours.

We next examine whether these same results carry over to the general curing case, or whether there are any notable differences. We should note that the curing and initialization policies are by nature not identical. They do, however, mirror one another in such a way that we can at least provide a surface-level comparison between the two setups.

\subsection{Delta-Curing Trials}\label{sec:curing_trials}

We consider the one-sided curing problem of trying to find a solution for~\eqref{eq:obj_curing}, for some time $n$. As in the initialization case, we compare the effectiveness of the various curing policies presented in Section~\ref{sec:heuristics}. A brief overview of these policies is found in Table~\ref{tab:curing_strats}.

\begin{table}[t!]
	\centering
	\caption{Delta-Curing Strategies}\label{tab:curing_strats}
	\begin{tabular}{cc}\hline\hline	
		$ \text{(i)} $ & Constrained gradient descent on a simplex: \\
			& Find $ \Delta_{b,i}(n) $ using Algorithm~\ref{alg:grad_descent} on the function $E[\tS_n|\F_{n-1}]$ \\
\hline
		$ \text{(ii)} $ & Uniformly allocate the budget to all nodes: \\
			& $ \Delta_{b,i}(n) =\frac{\budgetB}{N}  $ \\
	\hline
		$ \text{(iii)} $ & Uniformly allocate the budget to inner nodes: \\
			& $ \Delta_{b,i}(n)=\begin{cases}
						\frac{1}{|V'|}\budgetB	&\text{ if $i\in V'$} \\
						0						&\text{ else} \\
			\end{cases} $, where $V' = \{k\in V|\N_k'\not\subset\N_j',\forall j\in V\}$ \\
	\hline
		$ \text{(iv)} $ & Ratio of degree, closeness centrality, and super urn proportion (inner): \\
			& $	\Delta_{b,i}(n)=\begin{cases}
					\frac{|\N_i|C_iS_{i,n-1}}{\sum_{j\in V'}|\N_j|C_jS_{j,n-1}}\budgetB	&\text{ if $i\in V'$} \\
					0												&\text{ else}
	\end{cases} $\\
	\hline
		$ \text{(v)} $ & Allocate as in $ \text{(iii)} $ using $V'$ from Algorithm~\ref{alg:node_target} \\
	\hline
		$ \text{(vi)} $ & Allocate as in $ \text{(iv)} $, using $V'$ from Algorithm~\ref{alg:node_target} \\
	\hline
		$ \text{(vii)} $ & Allocate as in $ \text{(iii)} $ using $V'$ from Algorithm~\ref{alg:node_target_dense} (without Optional)\\
	\hline
		$ \text{(viii)} $ & Allocate as in $ \text{(iv)} $, using $V'$ from Algorithm~\ref{alg:node_target_dense} (without Optional)\\
	\hline
		$ \text{(ix)} $ & Ratio of degree, closeness centrality and super urn proportion: \\
			& $ B_i=\frac{|\N_i|C_iS_{i,n-1}}{\sum_{j\in V}|\N_j|C_jS_{j,n-1}}\budgetB $ \\
\hline\hline
	\end{tabular}
\end{table}

\begin{figure}
\centering{
  	\subfigure[Results of curing trials for low density Barabasi-Albert network depicted in~\ref{fig:BANet}]{ \includegraphics[width=0.7\linewidth]{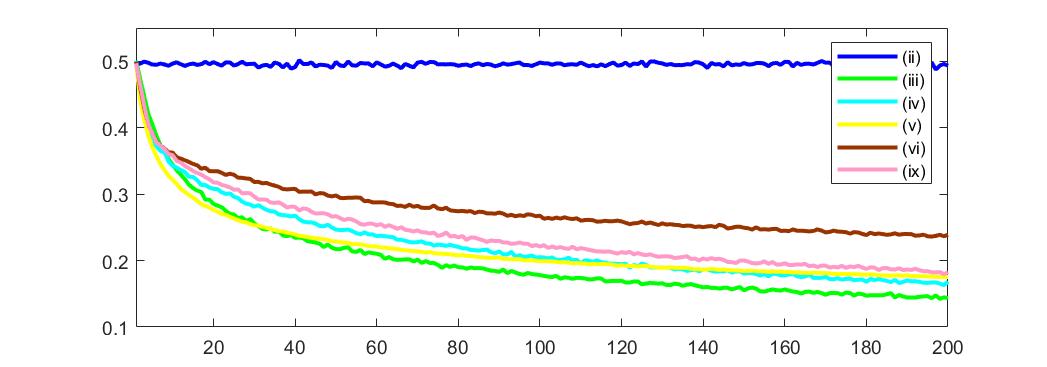}\label{fig:BACureTrials}}
}
\\
{
  	\subfigure[Results of curing trials for high density Barabasi-Albert network depicted in~\ref{fig:BADNet}]{ \includegraphics[width=0.7\linewidth]{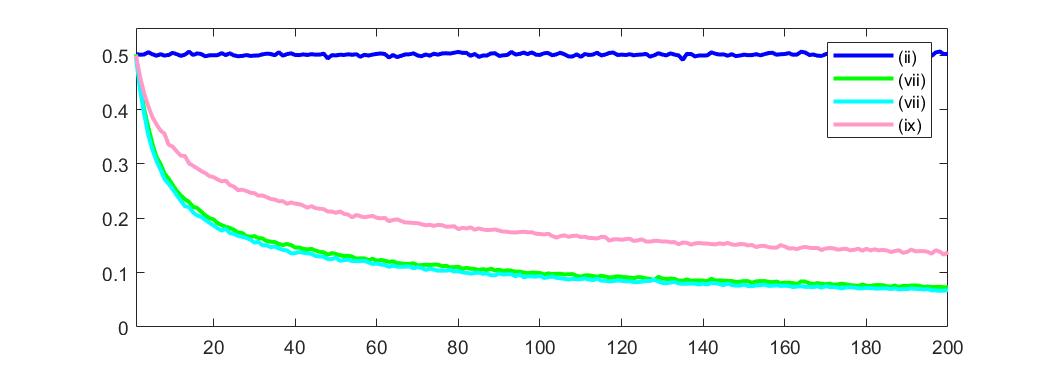}\label{fig:BADCureTrials}}
}
\\
{
  	\subfigure[Results of curing trials for Facebook network depicted in~\ref{fig:FaNet}]{ \includegraphics[width=0.7\linewidth]{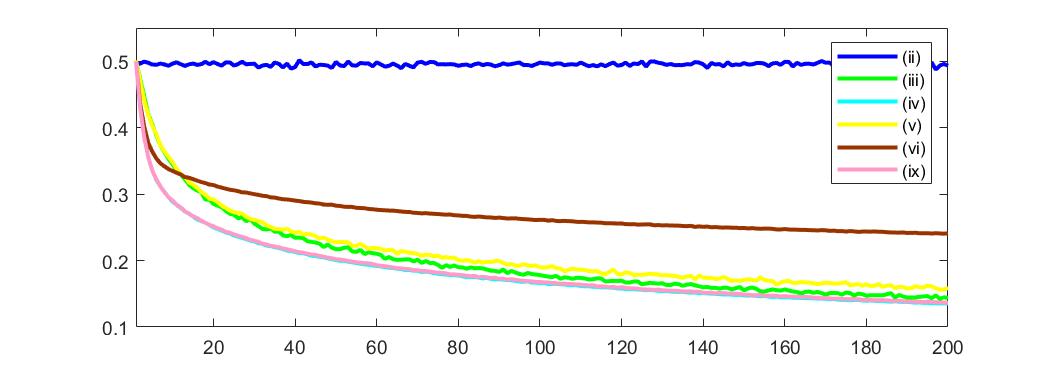}\label{fig:FaCureTrials}}
}
\caption[Comparison of curing policies]{Plots of empirical average infection rate, $\tI_n$, for various curing strategies across multiple networks. Lower values indicate lower levels of infection. In each case the network was initialized with $R_i=B_i=10$, $\forall i\in V$, and the Polya network contagion process was run over 500 trials with a fixed budget $\budgetR=\budgetB=10N$, with $\Delta_r=10$ fixed.}
\label{fig:CureTrials}
\end{figure}

We begin with a uniform initialization, where each node is allocated a fixed number of red balls, $R$, and black balls, $B$, with $R=B$. We then run the Polya network contagion process. Before each draw is performed, we assign a number of curing resources to each node $\Delta_{b,i}(n)$ in accordance with the desired policy from Table~\ref{tab:curing_strats}, while the infection resources are assigned uniformly, $\Delta_{r,i}(n)=\frac{\budgetR}{N}$ (with a fixed budget $\budgetB=\budgetR$ for each time $n$). As in the initialization case, we record the draw value, $Z_{i,n}$ for each node $i \in V$ at each time $n$. The draw values are averaged over the number of trials and the set of nodes to obtain a measure for the empirical average infection rate, $\tI_n$. The results of these simulations are depicted in Figure~\ref{fig:CureTrials}.

We first remark that unlike the initialization case, targeting a reduced subset of nodes based on centrality does not always lead to the best performance. This is evidenced by both Figures~\ref{fig:BACureTrials} and~\ref{fig:FaCureTrials}, where strategy $(vi)$ appears to perform worst of all the heuristics apart from the no control case. In fact, between the three network setups, the heuristics vary wildly in performance, with the best performing heuristic from~\cite{MH:17} performing extremely well in Figure~\ref{fig:FaCureTrials}, but performing significantly worse than the other techniques in~\ref{fig:BADCureTrials}. The various strategies also perform to different degrees of success for different values of $n$ (see strategy $(vi)$ in Figure~\ref{fig:FaCureTrials}).

Between the three networks, the strategies that seem to perform the most consistently well are strategies that target only the inner set of nodes (strategies $(iii)$, $(iv)$ and the analogous strategies for denser networks, $(vii)$ and $(viii)$). It is interesting that these strategies would perform most consistently for the Delta-curing problem, but would perform far from optimally for the one-time initialization policy as seen in Figure~\ref{fig:InitTrials}. The main difference between these setups is that the solutions to the Delta-curing problem are applied continually as the infection evolve, while solutions to the initialization problem are applied once. Obviously, some of the benefits of targeting a greatly reduced set of nodes, as in seen in Figure~\ref{fig:BANetAlg2}, are lost when they do not properly account for changes in the urn compositions. Interestingly, for denser networks, targeting a greatly reduced set of nodes as in Figure~\ref{fig:BADNetAlg3} was still very effective, most likely as a direct result of the high network density.

\begin{figure}
\centering{
  	 \includegraphics[width=0.3\linewidth]{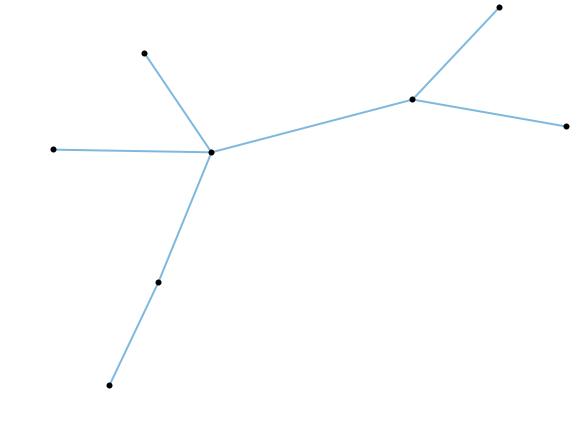}\label{fig:GradCompNet}
}
{
	\includegraphics[width=0.55\linewidth]{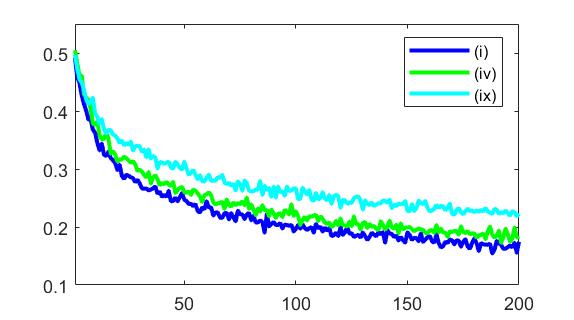}\label{fig:GradCompPlot}
}
\caption[Comparison of curing policies to gradient descent method]{Empirical average infection rate, $\tI_n$, for small 8 node test network. The performance of gradient descent over the expected network exposure, $E[\tS_n|\F_{n-1}]$, was compared to the two best performing heuristics from Figure~\ref{fig:BACureTrials}. The network was initialized with $R_i=B_i=10$, $\forall i\in V$, and the Polya network contagion process was run over 500 trials with a fixed budget $\budgetR=\budgetB=80$, with $\Delta_r = 10$ fixed.}
\label{fig:GradComp}
\end{figure}

One notable exclusion from Figure~\ref{fig:CureTrials} is the fact that none of the three subplots show the performance obtained using gradient descent, in part due to the highly complex nature of the expected network exposure $E[\tS_n|\F_{n-1}]$, which is the function we are minimizing over the curing parameters. While any large scale gradient-descent simulations are somewhat infeasible as a result, it is still possible to provide a comparison on a smaller scale, in order to see how far from `optimal' our various heuristics are.\footnote{As noted earlier, we say `optimal' since the gradient-descent algorithm only converges to an optimal solution for the one-step problem of minimizing $E[\tS_n|\F_{n-1}]$ over the curing resources.} 
We provide such a comparison in Figure~\ref{fig:GradComp}. While it is clear that the gradient descent does marginally improve the empirical average infection rate, we see that strategy $(iv)$ still performs comparably, at a fraction of the computational cost.


\section{Conclusion}\label{ch:conclusion}

We considered the (infinite-memory) Polya network contagion model and formulated (preemptive and reactive) control policies to mitigate the spread of infection. We framed these problems in two lights: both as a one-time initialization problem, and as a continual control problem that is an extension of contagion curing methods in prior work. We provided simplified metrics for the one-sided finite-horizon version of these problems and proved the existence of optimal policies. We developed effective heuristics that functioned primarily by modifying the scope of solutions to only a limited number of nodes. We then compared these heuristics to the one-step optimal policies that were obtained through use of gradient descent. We established underlying results for the initialization problem proving that optimal initialization strategies would only allocate resources to inner nodes. We showed how network symmetry impacts the allocation of resources for optimal policies. We also demonstrated that the proposed heuristics perform comparably to the one-step optimal policy in certain situations. Finally, we framed the Delta-curing problem as a two-person zero sum game on the expected network exposure, and we proved the existence of a Nash equilibrium. We borrowed ideas from the initialization problems in order to implement heuristics for the Delta-curing problem and demonstrated that these strategies can perform near-optimally at a fraction of the computational expense. Future directions include analogous optimization problems when underlying resource distribution is hidden from the observer, studying the initialization problem in a game-theoretic framework, and investigating the spread of the Polya contagion process in networks with unknown topologies.



\section*{Acknowledgement} The authors sincerely thank Tam\'{a}s Linder for considerably shortening the proof of Lemma~\ref{lma:315}.

\end{spacing}

\bibliographystyle{IEEEtran}
\bibliography{alias,GH-add,Main-add}


\end{document}